\input amstex
\input amsppt.sty
\magnification\magstep1
\input epsf
\def\ni\noindent
\def\sbs{\subset}

\def\as{\operatorname{asdim}}
\def\X{\operatorname{X-dim}}
\def\asdim{\operatorname{asdim}}
\def\diam{\operatorname{diam}}
\def\dim{\operatorname{dim}}

\def\R{\text{\bf R}}

\def\Q{\text{\bf Q}}
\def\Z{\text{\bf Z}}

\def\N{\text{\bf N}}

\def\sC{\Cal C}

\def\sF{\Cal F}

\def\sV{\Cal V}

\def\sU{\Cal U}
\def\Ucal{\Cal U}
\def\sW{\Cal W}

\def\p{\partial}

\hoffset= 0.0in
\voffset= 0.0in
\hsize=32pc
\vsize=38pc
\baselineskip=24pt
\NoBlackBoxes
\topmatter
\author
A. Dranishnikov
\endauthor

\title
Cohomological approach to asymptotic dimension
\endtitle
\abstract We introduce the notion of asymptotic cohomology based
on the bounded cohomology and define cohomological asymptotic
dimension $\as_{\Z} X$ of metric spaces. We show that it agrees
with the asymptotic dimension $\as X$ when the later is finite.
Then we use this fact to construct an example of a metric space
$X$ of bounded geometry with finite asymptotic dimension for which
$\as(X\times\R)=\as X$. In particular, it follows for this example
that the coarse asymptotic dimension defined by means of Roe's
coarse cohomology is strictly less than its asymptotic dimension.
\endabstract

\thanks The author was partially supported by NSF grant DMS-061519
\endthanks

\address University of Florida, Department of Mathematics, P.O.~Box~118105,
358 Little Hall, Gainesville, FL 32611-8105, USA
\endaddress


\email  dranish\@math.ufl.edu
\endemail

\endtopmatter

\document

\head \S0 Introduction\endhead

\

Gromov proposed to study discrete groups as large scale geometric
objects. He introduced several asymptotic invariants of finitely
generated groups [Gr1]. Among them there is the notion of
asymptotic dimension which proved to be useful for the Novikov-type
conjectures [Yu],[Ba],[CG],[Dr2],[DFW]. The asymptotic
dimension differs from any other known dimension invariant of
discrete groups. Thus, it differs from the geometric dimension
$gd(\Gamma)$ [Br] for every finite group $\Gamma$ since the geometric dimension
is infinite for groups with torsions and the asymptotic dimension
of finite groups is zero. If
one takes into account only torsion free groups then the
distinction between these two dimensions is not obvious. For almost
all known torsion free groups where both invariants are computed
there is the equality $\as\Gamma=gd(\Gamma)$. An exception is
Gromov's example of a group $\Gamma_0$ with finite classifying
space $B\Gamma_0$ that contains an expander (in some weakly coarse
sense) [Gr2]. This group does not admit a coarse embedding into
a Hilbert space and hence [HR], $\as\Gamma_0=\infty$
whereas $gd(\Gamma_0)<\infty$.

For discrete groups it is known that modulo some potential
discrepancy in the case of geometric dimension three
(Eilenberg-Ganea problem) there is the equality
$gd(\Gamma)=cd(\Gamma)$ where $cd(\Gamma)$ is the cohomological
dimension of $\Gamma$ with integral coefficients [Br]. This equality
gives an effective way for computation of $gd(\Gamma)$. A similar
situation happens to be in the classical dimension theory of locally compact
metric spaces. One of the basic facts there is the
Alexandroff Theorem which states that
the covering dimension agrees with the integral cohomological
dimension: $\dim X=\dim_{\Z}X$, provided $\dim X<\infty$.
The cohomological approach in the classical dimension theory
very often allows one to compute the covering
dimension or to reduce the computation to a calculation of
cohomological dimensions of $X$ with respect to the coefficient groups
$\Q$, $\Z_p=\Z/p\Z$, and $\Z_{p^{\infty}}=\lim_{\rightarrow}\Z_{p^k}$ for
prime $p$. Also it gives an exact formula for dimension of the product
(Bockstein's formula [Dr3]). In particular the cohomological approach
to dimension allows to prove Morita's theorem which states that
$\dim(X\times\R)=\dim X+1$ for general topological spaces $X$.

In the light of all this it would be useful to have a cohomological
approach to the asymptotic dimension. The right approach would give
answer to many naive still open questions in asymptotic dimension
theory, like, {\it Does the asymptotic Morita theorem hold
true: $\as(\Gamma\times\Z)=\as\Gamma+1$ for finitely generated groups
$\Gamma$?}

Perhaps the first attempt to define asymptotic dimension
cohomologically was made in [Dr1] by means of Roe's coarse
cohomology theory $HX^*$. We redefine the corresponding dimension
in this paper (\S4) and call it the coarse cohomological dimension
${\X}_G$. It turns out that the coarse cohomological dimension
${\X}_G$ does not always agree with the asymptotic dimension
even for proper metric spaces of bounded geometry and when the
later is finite (\S5). Thus, an asymptotic analog of Alexandroff
Theorem does not hold for this dimension. We recall that the most
elegant argument for the Alexandroff theorem uses the
representation theorem for the \v Cech cohomology. In view of our
result the macro-micro analogy suggests that Roe's coarse
cohomology theory lacks Brown's representability.

Moreover, it turns out that the asymptotic Morita theorem does not
hold true for general proper metric spaces of bounded geometry. We
construct a counterexample in \S5 based on a new cohomological
approach to asymptotic dimension and an idea from elementary
number theory. In \S6 we show how this example can be turned into
a large scale fractal-like space. We present there a general
construction of large simplicial complexes of prescribed shapes on
all scales.

A new cohomological approach to the asymptotic dimension is based
on cohomology groups defined by bounded cochains. Using bounded
cohomology we define an asymptotic cohomological dimension
$\as_{\Z}$ with integral coefficients in \S2 and we show that
$\as_{\Z}X=\as X$ for proper metric spaces of bounded geometry
with $\as X<\infty$ in \S3. So this approach has an asymptotic
analog of the Alexandroff Theorem and it also gives a hope for
developing an analog of the Bockstein theory.

For finite coefficients $F$ the bounded cohomologies coincide with
the standard cohomologies. We use this to show (\S4) that
$\as_FY=X{-}\dim_FY$ in this case.  In view of this equality, our approach to
the asymptotic dimension still gives a hope for the formula
$\as\Gamma=cd(\Gamma)$ for torsion-free finitely presented groups with
finite asymptotic dimension. Now in order to derive this formula it
suffices to show that in the case of a finitely presented group
$\Gamma$, $\as_{\Z}\Gamma=\as_{\Z_p}\Gamma$ for some $p$. We recall
that in the world of compacta this holds true for all
locally nice spaces [Dr3]. It is still unclear whether finitely presented
groups are sufficiently nice among all discrete proper metric spaces.

\head \S1 Preliminaries\endhead

\

{\bf Asymptotic dimension.} Gromov defined the asymptotic
dimension of a metric space $X$ as follows [Gr1].

DEFINITION.  {\it A metric space $X$ has asymptotic dimension $\leq n$ if,
for every $d>0$, there is an $R$ and $n+1$ $d-$disjoint, $R-$bounded families
$\Ucal_0, \Ucal_1, \ldots, \Ucal_n$ of subsets of $X$ such that
$\cup_{i=0}^n \Ucal_i$ is a cover of $X$.}

We say that a family $\Ucal$ of subsets of $X$ is $R-$bounded if
$\sup \{ \diam{U} | U \in \Ucal \} \leq R$.  Also, $\Ucal$ is
said to be $d-$disjoint
if $d(x,y) > d$ whenever $x\in U$, $y \in V$, $U \in \Ucal$,
$V \in \Ucal$, and $U \neq V$.

It is known that $\as X\le n$ if and only if $X$ admits a uniformly bounded
open cover $\sU$ of multiplicity $\le n+1$ with an arbitrary large
Lebesgue number $L(\sU)=\inf_{x\in X}\sup_{U\in\sU}d(x,X\setminus U)$.

The notion of asymptotic dimension is a coarse invariant (see \cite{Ro1}
\cite{Ro2}). Therefore the invariant $\as\Gamma$ is well-defined
for every finitely generated group (in fact for every countable
group $\Gamma$ \cite{DS}).

A metric space $X$ is called {\it proper} if every closed ball
$B_r(x)\subset X$ is compact.

\

We recall that a cover $\sV$ of a space $X$ is called a {\it
refinement} of a cover $\sU$, (the notation for this relation is
$\sV\prec\sU$), if for every $V\in\sV$ there is $U\in\sU$ such
that $V\subset U$. Note that if the mesh of $\sV$ is less than the
Lebesgue number $L(\sU)$ of $\sU$, then $\sV\prec\sU$. A
refinement defines a simplicial map of the nerves $q:N(\sV)\to
N(\sU)$ which is called a {\it refinement map}.

Let $\sC$ be a collection of sets in $X$ and let $A\subset X$.
The star of $A$ with respect to $\sC$ is the set
$$
St(A,\sC)=\bigcup_{C\in\sC,C\cap
A\ne\emptyset}C.$$
If $\tau$ is the collection of simplices of a simplicial complex $K$,
then for every vertex $v\in K$, $St(v,\tau)=St(v,K)$ is the combinatorial
star of the vertex. The open star of a vertex $Ost(v,K)= Int(St(v,K))$
can be also defined as the complement to the link $St(v,K)\setminus Lk(v,K)$.

A cover $\sV$ of a space $X$ is called a {\it star refinement}
of a cover $\sU$, $\sV\prec\prec\sU$, if for every $V\in\sV$ there
is $U\in\sU$ such that $St(V,\sV)\subset U$. The corresponding simplicial map
$q:N(\sV)\to N(\sU)$ is called a {\it star refinement map}.

{\bf Anti-\v{C}ech approximation.} A sequence of uniformly bounded
locally finite open coverings $\{\sU_i\}$ of a metric space $X$ is
called an {\it anti-\v{C}ech approximation} [Ro1] for $X$ if
$\sU_i\prec\sU_{i+1}$ for all $i$ and the Lebesgue number
$L(\sU_i)$ tends to infinity. Let $N_i=N(\sU_i)$ denote the nerve
of $\sU_i$. Then every anti-\v{C}ech approximation defines a
direct system of locally finite simplicial complexes with
refinement maps as the bonding maps:
$$
N_1 @>q^1_2>> N_2 @>q^2_3>> \dots @>>> N_k @>q^k_{k+1}>> N_{k+1}
@>>>\dots\ \ .
$$
For a locally finite open cover $\sU$ of a metric space $X$ let
$p:X\to N(\sU)$ denote a projection to the nerve defined by the
partition of unity $\{\phi_U\}$ with $\phi^{-1}_U(0)=X\setminus
U$. We call such projection {\it canonical} if it is defined by
the following partition of unity:
$$\phi_U(x)=\frac{d(x,X\setminus U)}{\sum_Vd(x,X\setminus V)}.$$
We recall that a cover $\sU$ is called {\it irreducible} if for
every $U\in\sU$ the family $\sU\setminus\{U\}$ is not a cover. If
the cover $\sU$ is irreducible then all vertices of $N(\sU)$ are in
the image of $p$. For a map $p$ with such property we say that it is
{\it essentially surjective}. We always will assume that the covers
$\sU_i$ in the definition of an anti-\v Cech approximation are
irreducible.

A {\it uniform} metric on a simplicial complex $K$ is the metric
restricted from the Hilbert space $\ell_2(K^{(0)})$ under the
natural realization of $K$. The geodesic metric on $K$ induced
from the uniform metric is called {\it uniform geodesic}. Usually
we normalize this metric by $1/\sqrt{2}$ to have the length of
every edge in all simplices equal one.

A map $f:X\to K$ of a metric space to a simplicial complex is
called {\it uniformly cobounded} if there is $D>0$ such that
$diam(f^{-1}(\Delta))\le D$ for every simplex $\Delta$.

Let $f:X\to Y$ be a map between metric spaces. We denote the
number (if it exists)
$$Lip(f)=\sup\{\frac{d_Y(f(x),f(y))}{d_X(x,y)}\mid x,y\in X, x\ne y\}$$ and
call it {\it the Lipschitz constant} of $f$. Every number
$\lambda\ge Lip(f)$ is called {\it a Lipschitz constant} for $f$.

We note that for every irreducible uniformly bounded locally
finite open cover $\sU$ of a metric space $X$ with the Lebesgue
number $L$ the canonical projection $p:X\to N(\sU)$ to the nerve
is $c_n/L$-Lipschitz and uniformly cobounded where $n$ is the
multiplicity of $\sU$ [BD1]. Moreover, it is a quasi-isometry
provided $X$ is geodesic and the nerve is taken with geodesic
metric:

We recall that $p:X\to Y$ is a {\it $(l,D)$-quasi-isometric
embedding} if there are constants $\lambda\ge 1$ and $D$ such that
$$
\frac{1}{l}d_X(x,y)- D\le d_Y(p(x),p(y))\le l d_X(x,y)+D.
$$
An  $(l,D)$-quasi-isometric embedding $p:X\to Y$ is called a {\it
quasi-isometry} if there is a constant $r$ such that the image
$p(X)$ is $r$-dense in $Y$, i.e. $N_{r}(p(X))=Y$. We use the
notations $N_r(A)=\{x\in X\mid d(x,A)\le r\}$ and $ON_r(A)=\{x\in
Y\mid d(x,A)<r\}$, $r>0$ for closed and open $r$-neighborhoods of
the set $A$ in a metric space $Y$. Also we consider "negative"
neighborhoods $N_{-r}(A)=Y\setminus N_r(X\setminus A)$, $r>0$.

 The following lemma can be extracted from [BD1] and {Dr4].

\proclaim{Lemma 1.1} Given $n$ there is a constant $C_n$ such that
for every uniformly bounded cover $\sU$ of a proper geodesic
metric space $X$ with multiplicity of $\sU\le n+1$ the canonical
projection $p:X\to N(\sU)_{UG}$ to the nerve is a
$\epsilon$-Lipschitz $(l,3)$-quasi-isometry with $r=1$
where $\epsilon=C_n/L(\sU)$ and $l=C_nmesh(\sU)$.
\endproclaim

Note that if the multiplicity of coverings $\sU_i$ in an anti-\v
Cech approximation of $X$ is bounded from above then the canonical
projections $p_i:X\to N(\sU_i)$ are $\epsilon_i$-Lipschitz with
$\lim\epsilon_i=0$. Thus, a metric space $X$ has $\as X\le n$ if
and only if it admits an anti-\v Cech approximation
$\sU_1\prec\sU_2\prec\sU_3\prec\dots$ with the multiplicity of
$\sU_i$ bounded from above by $n+1$. In terms of projections to
the nerves it equals to the property that $X$ admits a sequence
$p_i:X\to N_i$ of uniformly cobounded $\epsilon_i$-Lipschitz maps
to uniform $n$-dimensional complexes with $\lim\epsilon_i=0$.

We recall that a metric space $X$ is said to be {\it of bounded
geometry} (on large scale) if for every $R$ the 1-capacity of
$R$-balls in $X$ is uniformly bounded from above. It was shown in
[Ro1] that given $\lambda>0$, every metric space of bounded
geometry $X$ admits an open cover of finite multiplicity with the
Lebesgue number greater than $\lambda$. Thus every metric space of
bounded geometry admits an anti-\v Cech approximation.

Let $p:K\to N$ be map between simplicial complexes. A simplicial
map $q:K\to N$ is called a {\it simplicial approximation} of $p$
if $p^{-1}(\Delta)\subset q^{-1}(\Delta)$ for every simplex
$\Delta\subset N$. This is equivalent to the condition
$p^{-1}(L)\subset q^{-1}(L)$ for every subcomplex $L\subset N$.

\proclaim{Proposition 1.2} Every proper geodesic metric space $X$
with $\as X\le n$ admits an anti-\v{C}ech approximation $\{\sU_i,
q^i_{i+1}\}$ with $n$-dimensional locally finite nerves $N_i$ and
essentially surjective projections $p_i:X\to N_i$ such that
\roster \item{} there are bonding maps $p^i_{i+1}:N_i\to N_{i+1}$
with $p_{i+1}=p^i_{i+1}\circ p_i$ for all $i$, \item{}
$(p_{i+1}^i)^{-1}(K)$ is a subcomplex for every subcomplex
$K\subset N_{i+1}$ , \item{} simplicial maps $q^i_{i+1}:N_i\to
N_{i+1}$ are simplicial approximations of $p^i_{i+1}$, \item{}
$Lip(p^i_{i+1})<1/2$, and \item{} $\sU_i\prec\prec\sU_{i+1}$ for
all $i$ and $q^i_{i+1}$ is a star refinement map.
\endroster
\endproclaim
\demo{Proof} We construct these coverings and maps by induction.
Assume that a sequence of covers
$\sU_1\prec\prec\sU_2\prec\prec\dots\sU_k$ together with the maps
$p_i$, $p^i_{i+1}$, and $q^i_{i+1}$ satisfying conditions (1)-(5)
is constructed.
We assume that $p_1$ and all $p^i_{i+1}$ are canonical projections
to the nerves. By Proposition 1.1 $N_k$ is quasi-isometric to $X$
and hence $\as N_k\le n$. We consider a uniformly bounded cover
$\sV'$ of $N_k$ of multiplicity $n+1$ with the Lebesgue number
$\ge 2C_n+2\ge 5$. Then $N_{-1}(\sV')=\{N_{-1}(V')\mid V'\in\sV'\}$
has the Lebesgue number $\ge 2C_n+3$. We define
$\sV=\{Int(St(V,N_k))\mid V\in N_{-1}(\sV')\}$.
Then $\sV\prec\sV'$ and hence, the nerve of $\sV$ is $n$-dimensional
and $\{Ost(v,N_k)\mid v\in N_k^{(0)}\}\prec\prec\sV$. Define
$\sU_{k+1}=p^{-1}_k\sV$. Then $$\sU_k=p^{-1}_k\{Ost(v,N_k)\mid v\in
N_k^{(0)}\}\prec\prec p^{-1}_k\sV=\sU_{k+1}.$$
Let
$N_{k+1}=N(\sV)=N(\sU_{k+1})$, let $p^k_{k+1}:N_k\to N_{k+1}$ be
the canonical projection to the nerve, and let
$p_{k+1}=p^k_{k+1}p_k$. We define
$q^k_{k+1}(v)$ to be an element $Int(St(V,N_k))$ of $\sV$, $V\in\sV'$ such that
$St(v)\subset V$. By Proposition 1.2 $p^k_{k+1}$ is
$C_n/L(\sV)$-Lipschitz. Thus, $Lip(p^k_{k+1})\le 1/2$ by the
choice of $\sV$. Thus the conditions (1),(4), and (5) are
satisfied.

To verify (2) it suffices to show that the preimage
$(p^k_{k+1})^{-1}(\Delta)$ is a subcomplex in $N_k$ for every
simplex $\Delta$ in $N(\sV)$.  Note that
$$(p^k_{k+1})^{-1}(\Delta)=N_k\setminus\bigcup_{W\notin\Delta^{(0)}}W=
N_k\setminus\bigcup_W Int(K_W)=\bigcap_W(N_k\setminus Int K_W)$$
where $K_W=St(V_W,N_k)$ is a subcomplex of $N_k$. Hence
$(p^k_{k+1})^{-1}(\Delta)$ is a subcomplex of $N_k$ as the
intersection of subcomplexes.

Finally we check (3). Let $p^k_{k+1}(z)\in\Delta=[V_1,\dots,
V_j]$. Then if $z\in V\in\sV$, it follows that $V$ equals one of
$V_i$, $i\le j$. Let $z\in[v_1,\dots, v_s]\subset N_k$. By the
definition $z\in St(v_i)\subset q^k_{k+1}(v_i)$. Therefore,
$q^k_{k+1}(v_i)\in\Delta^{(0)}$. Hence $q^k_{k+1}(z)\in\Delta$.
\qed
\enddemo

We will refer to such anti-\v Cech approximation as to {\it
regular} and will denote it as $$\{p_i:X\to N_i,q^i_{i+1},
p^i_{i+1}\}$$ where $p_i$ are $\epsilon_i$-Lipschitz with
$\lim\epsilon_i=0$. In some instances we will not assume the
condition (5). Note that the coverings $\sU_i$ can be recovered
from this data as $\sU_i=p_i^{-1}(Ost(v,N_i))$ where $Ost(v,N)$ is
the open star of a vertex $v$ in a simplicial complex $N$.

\head \S2 Definition of $\as_{\Z}$\endhead

\

{\bf Bounded cohomology.} Let $K$ be a simplicial complex. An
integral cochain $\phi:C_m(K)\to \Z$ is bounded if there is a
constant $b<\infty$ such that $|\phi(\sigma)|<b$ for all
$m$-simplices $\sigma$ in $K$. Clearly, the coboundary
$\delta\phi$ of a bounded cochain is a bounded cochain. Cohomology
groups defined by means of bounded cochains $C_b^*$ are called
{\it bounded cohomologies} of $K$ and denoted as $H_b^*(K)$.
Clearly, the inclusion $C_b^*\to C^*$ is a chain map. Therefore
there is a natural homomorphism $H_b^*(K)\to H^*(K)$. Every
simplicial map $\phi:K\to L$ induces a homomorphism $H_b^*(L)\to
H_b^*(K)$. This definition can be extended to any coefficient
group with a norm such as $\Q$ or $\R$ and their subgroups. For a
subcomplex $L\subset K$ one can define a relative bounded
cohomology groups $H_b^*(K,L)$ by considering the relative
cochains. We note that for bounded cohomology there are the exact
sequence of pair
$$
\dots @<<< H^i_b(L) @<<< H^i_b(K) @<<< H^i_b(K,L) @<<<
H^{i-1}_b(L) @<<<\dots,
$$
exact sequence of triple
$$
\dots @<<< H^i_b(L,A) @<<< H^i_b(K,A) @<<< H^i_b(K,L) @<<<
H^{i-1}_b(L,A) @<<<\dots
$$
where $A\subset L$ is a subcomplex, and the excision isomorphism
$$
H^i_b(K,B) @>\cong>> H^i_b(A,A\cap B)
$$
where $A,B\subset K$ are subcomplexes such that $K=A\cup B$.

{\bf Approximation by asymptotic polyhedra.} Let $X$ be a metric
space. We consider locally finite covers $\sU$ of $X$ by bounded
open sets such that the Lebesgue number $L_x(\sU)$ tends to
infinity as $x\to\infty$.  We say that $\sV$ is a refinement of
$\sU$ at infinity if there is $R>0$ such that for every $V\in\sV$,
$d(V,x_0)>R$, there is $U\in\sU$ with $V\subset U$. A refinement
at infinity $\sV \prec_{\infty} \sU$ defines a simplicial map
$\phi:K_{\sV}\to K_{\sU}$ between corresponding subcomplexes of
the nerves. All such covers with the relation $\prec_{\infty}$
form a directed set $aCov(X)$.

The following proposition is proven in [Dr1].
\proclaim{Proposition 2.1} If a metric space $X$ has $\as X\le n$
then the family $aCov(X)$ has a cofinal subfamily $aCov_n(X)$ that
consists of covers with $n$-dimensional nerves.
\endproclaim
A countable simplicial complex $K$ with a metric $d$ such that
every simplex is isometric to an affine simplex in a Hilbert space
is called  an {\it asymptotic polyhedron} if
$
\lim_{i\to\infty} Lip(\phi_{i})=0
$
where $K$ is the union of simplices $K=\cup_i\sigma_i$ and
$\phi_i:\sigma_i\to\Delta^{\dim\sigma_i}$ is the affine map to the
standard simplex.

We note that if a cover $\sU$ has bounded multiplicity then its
nerve $N$ admits a metric of an asymptotic polyhedron such that
the projection $p_{\sU}:X\to N$ is 1-Lipschitz (see [Dr1]
[DFW]).

{\bf Asymptotic bounded cohomology.} Let $K$ be a simplicial
complex. Denote by $\sF(K)$ the set of all finite subcomplexes of
$K$. We define the {\it bounded cohomology at infinity}
$AH_b^i(K)$ of a complex $K$ in the dimension $i$  as the direct
limit
$$
AH^i_b(K)=\lim_{\rightarrow}\{H_b^i(K\setminus F)\mid F\subset K,
F\in\sF(K)\}.
$$
Let $X$ be a metric space and let $\sV\prec_{\infty}\sU$ be covers
as above. Then there is a well-defined homomorphism
$AH_b^*(N(\sU))\to AH_b^*(N(\sV))$ between the bounded cohomology
at infinity of nerves.

For a metric space $X$ we define its  {\it asymptotic bounded
cohomology} as the direct limit
$$
AH^i(X)=\lim_{\rightarrow}\{AH^i_b(N(\sU))\mid \sU\in aCov(X)\}.
$$
Thus, it can be defined as
$$
AH^i(X)=\lim_{\rightarrow}\{H_b^i(N(\sU)\setminus F)\mid\sU\in aCov(X), F\in\sF(N(\sU))\}.
$$
This definition can be extended to any coefficient group with the
norm.

Let $L$ be a subcomplex of a simplicial complex $K$.  We define
the bounded cohomology at infinity $AH_b^*(K,L)$ of the pair
$(K,L)$ in the dimension $i$ as the direct limit
$$
AH_b^i(K,L)=\lim_{\rightarrow}\{H_b^i(K\setminus F,L\setminus
F)\mid F\subset K, F\in\sF\}.
$$
Let $X$ be a metric space and let $Y\subset X$ be a subset. For a
cover $\sU$ of $X$ we denote by $N(\sU|_Y)$ the nerve of the cover
$\sU$ restricted to $Y$, $\sU|_Y=\{U\cap Y\mid U\in\sU\}$. Clearly,
$N(\sU|_Y)\subset N(\sU)$ and a refinement $\sV\prec\sU$ defines a
simplicial map of pairs $(N(\sV),N(\sV|_Y))\to(N(\sV),N(\sU|_Y))$.
Then we define a relative asymptotic cohomology as the direct limit
$$
AH^i(X,Y)=\lim_{\rightarrow}\{AH_b^i(N(\sU),N(\sU|_Y))\mid \sU\in
aCov(X)\}.
$$

Let $Y\subset X$ be a subset of a metric space. For a cover $\sU$ of
$X$ we denote by $\sU_Y=\{U\in\sU\mid U\cap Y\ne\emptyset\}$. By
$aCov_X(Y)$ we denote the set of locally finite covers $\sV$ of $Y$
by bounded open sets from $X$ with $lim_{x\to\infty,x\in
Y}L_x(\sV)=\infty$. Thus, for every $\sU\in aCov(X)$ we have
$\sU_Y\in aCov_X(Y)$. Note that every cover  $\sV\in aCov_X(Y)$ can
be enlarge to a cover $\sU\in aCov(X)$, $\sV\subset\sU$, in such a
way that $\sU|_Y=\sV|_V$.

\proclaim{Proposition 2.2} The family $\{\sU|_Y\mid \sU\in
aCov(X)\}=\{\sU|_Y\mid \sU\in aCov_X(Y)\}$ is cofinal in $aCov(Y)$
where $Y\subset X$ is taken with the restriction metric.
\endproclaim
\demo{Proof} Let $\sW\in aCov(Y)$. For every $W\in\sW$ we define $
\tilde W=\cup_{x\in W}B_{r_x/3}(x)$ where $r_x$ is the supremum of
$r$ such that $B_r(x)\cap Y\subset W$ and $B_r(x)$ is the $r$-ball
in $X$. Clearly, $\tilde W\cap Y=W$. Let $\tilde\sW=\{\tilde W\mid
W\in\sW\}$. We note the Lebesgue number of $\tilde\sW$ restricted to
 $Y$ tends to infinity. Thus, $\tilde W\in aCov_X(Y)$
\qed
\enddemo
REMARK. Note that $N(\tilde\sW)=N(\sW)$. Since $\{\tilde \sW\mid
\sW\in aCov(Y)\}$ is cofinal in $aCov_X(Y)$, we obtain
$$AH^*(Y)=\lim_{\rightarrow}\{AH_b^*(N(\sU|_Y));\sU\in aCov(X)\}=
\lim_{\rightarrow}\{AH_b^*(N(\sU_Y));\sU\in aCov(Y)\}.$$ Thus, the
asymptotic cohomology of a pair can be defined as
$$
AH^i(X,Y)=\lim_{\rightarrow}\{AH_b^i(N(\sU),N(\sU_Y))\mid \sU\in
aCov(X)\}.
$$

Since bounded cohomologies at infinity posses the exactness, the
exactness is preserved by direct limits, and in view of Proposition
2.2, there are exact sequences for asymptotic cohomology of pair
(and triple).

A triad $X,A,B$, $A\cup B=X$, is called {\it excisable} if the
family $ \{\sU_A\cap\sU_B\ \mid \sU\in aCov(X)\}$ is cofinal in
$aCov_X(A\cap B)$.

\proclaim{Proposition 2.3} For every excisable triad $X,A,B$ there
is the Mayer-Vietoris exact sequence
$$
\dots\to AH^n(A)\oplus AH^n(B)\to AH^n(X)\to AH^{n+1}(A\cap
B)\to\dots
$$
\endproclaim
\demo{Proof} We note that the Mayer-Vietoris sequence holds for
bounded cohomology and hence for bounded cohomology at infinity for
simplicial complexes. We consider the direct limit of these
Mayer-Vietoris sequences for the nerves $N(\sU)$, $\sU\in aCov(X)$:
$$
\to AH^n_b(N(\sU_A))\oplus AH^n_b(N(\sU_B))\to AH^n_b(N(\sU))\to
AH^{n+1}_b(N(\sU_A)\cap N(\sU_B))\to\ .
$$
Note that there are the inclusions $ N(\sU_A)\cap N(\sU_B)=
N(\sU_A\cap\sU_B).$ Since $ \{\sU_A\cap\sU_B\mid \sU\in aCov(X)\}$
is cofinal in $aCov_X(A\cap B)$, one can argue that
$$\lim_{\rightarrow}AH^*_b((N(\sU|_A)\cap N(\sU|_B))=
\lim_{\rightarrow}AH^*_b(N(\sU|_{A\cap B}))=AH^*(A\cap B).$$ \qed
\enddemo
There is a version of the Mayer-Vietoris sequence for pairs.

\proclaim{Proposition 2.4} For every excisable triad $A\cup B,A,B$
of subsets in a metric space $X$ there is the Mayer-Vietoris exact
sequence
$$
\dots\to AH^n(X,A)\oplus AH^n(X,B)\to AH^n(X,A\cap B)\to
AH^{n+1}(X,A\cup B)\to\dots\ .
$$
\endproclaim
A metric space $X$ is called {\it uniformly path connected} if there
is a monotone tending to infinity continuous function
$S:\R_+\to\R_+$, $S(0)=0$, such that every two points $x,x'\in X$
can be joined by a path $J$ with $diam(J)\le S(d(x,x'))$. We note
that every geodesic metric space is uniformly path connected.

\proclaim{Proposition 2.5} Let $W\subset X$ be an open subset of a
uniformly path connected metric space $X$. Then the triple $X,\bar
W,X\setminus W$ is excisable.
\endproclaim
\demo{Proof} Let $\sV\in aCov_X\p W$. Let $\sU'\in aCov(X)$ be an
enlargement of $\sV$, $\sV\subset\sU'$, and $\sU'|_Y=\sV|_Y$. Since
$X$ is uniformly path connected, the cover $\sU$ that consists of
components of sets from $\sU'$ belongs to $aCov(X)$. Then for every
$U\in\sU$ with $U\cap\bar W\ne\emptyset$ and $U\cap(X\setminus
W)\ne\emptyset$, we obtain $U\cap\p W\ne\emptyset$. Thus, $U$ is a
component of $U'$ with $U'\in\sV$. We have checked that $\sU|_{\bar
W}\cap \sU|_{X\setminus W}\prec\sV$. \qed
\enddemo

\proclaim{Proposition 2.6} Let $W\subset X$ be an open subset of a
uniformly path connected metric space $X$. Then there is the
excision isomorphism
$$
AH^n(X,X\setminus W)=AH^n(\bar W,\p W).
$$
\endproclaim
\demo{Proof} Let $\sU\in aCov(X)$. Note that there is an inclusion
$\p N(\sU_{\bar W})\subset N(\sU_{\bar W}\cap\sU_{X\setminus W})$
which is essentially onto. Moreover, if $\sV\prec\prec\sU$, then
there is a map $\xi:N(\sV_{\bar W}\cap\sV_{X\setminus W})\to\p
N(\sU_{\bar W})$ that make the diagram generated by $\xi$ and the
refinement map commutative. This implies that
$$\lim_{\rightarrow}AH_b^*(\p N(\sU_{\bar
W}))\to\lim_{\rightarrow}AH_b^*(N(\sU_{\bar W}\cap\sU_{X\setminus
W}))$$ is an isomorphism. In view of Proposition 2.5 the later limit
is equal to $\lim_{\rightarrow}AH_b^*(N(\sU_{\p W}))$. Therefore,
$$\lim_{\rightarrow}AH_b^*(N(\sU_{\bar W}),\p N(\sU_{\bar
W}))=\lim_{\rightarrow}AH_b^*(N(\sU_{\bar W}),N(\sU_{\p
W}))=AH^*(\bar W,\p W).$$ By the excision for bounded cohomology at
infinity of simplicial complexes we obtain
$$AH^*(X,X\setminus
W)=\lim_{\rightarrow}AH^*_b(N(\sU),N(\sU_{X\setminus W})=
\lim_{\rightarrow}AH^*_b(N(\sU_{\bar W}),\p N(\sU_{\bar W})).$$ \qed
\enddemo

 {\bf Asymptotic cohomological
dimension.} Let $X$ be a metric space we define its asymptotic
integral cohomological dimension as follows:
$$
\as_{\Z}X=\max\{n\mid AH^n(A,B)\ne 0\mid B\subset A\subset X\}.
$$
This definition can be extended to any coefficient group $G$ with a
semi-norm. The notation is $\as_GX$ when the semi-norm on $G$ is
specified.

REMARK. One can show that like in the case of cohomological
dimension in topology it suffices to consider $A=X$ in the above
definition. Indeed, from exact sequence of triple $B\subset
A\subset X$
$$
AH^{n+1}(X,A) @<<< AH^n(A,B) @<<< AH^n(X,B) @<<< AH^n(X,A) @<<< $$
it follows that if $AH^n(A,B)\ne 0$ for $n=\as_{\Z}X$, then
$AH^{n+1}(X,A)=0$ and hence $AH^n(X,B)\ne 0$.

\

\head \S3 Connection between $\as_{\Z}$ and $\as$\endhead

In this section we show that $\asdim_{\Z}$ agrees with $\asdim$ for
proper geodesic metric spaces provided the later is finite.

\proclaim{Proposition 3.1} Let $K\subset N$ be a subcomplex of a
uniform geodesic complex of dimension $n$. Suppose that $f:K\to L$
is the map to an $b$-bounded metric space such that
$Lip(f|_{\sigma})\le \lambda$, $\lambda\ge 1$, for all simplices
$\sigma\subset K$. Then $f$ is $c\lambda$-Lipschitz where $c$
depends on $n$ and $b$ only.
\endproclaim
\demo{Proof} Let $K$ be realized in the Hilbert space
$\ell_2(K^{(0)})$ spanned by the vertices of $K$. Then the
identity map $1_K:(K,d_N)\to(K,d_{\ell_2})$ is 1-Lipschitz. It
suffices to show that $f$ is $c\lambda$-Lipschitz with respect to
the Hilbert space metric. Note that the distance between two
disjoint simplices in $\Delta\subset \ell_2$ of dimension $\le n$
is at least $\sqrt{\frac{2}{n+1}}$. For every chain $0\le l\le
k\le k'$ we denote by $D_{l,k,k'}$ the union of standard simplices
$\Delta^k\cup\Delta^{k'}$ in $\ell_2$ with the intersection
$\Delta^k\cap\Delta^{k'}=\Delta^l$. Let $c_{l,k,k'}$ be the
Lipschitz constant of the identity map $id:(D_{l,k,k'},
d_{\ell_2})\to(D_{l,k,k'}, d_{geod})$ where $d_{geod}$ is the
intrinsic metric on $D_{l,k,k'}$ induced by the Euclidean metric.
Let $\bar c=\max\{c_{l,k,k'}\mid k'\le 2n+1\}$. We take
$c=b\sqrt{n+1}\bar c$. Then for every couple of points $x,x'\in K$
with $\|x-x'\|\ge\sqrt{\frac{2}{n+1}}$ we obtain
$d_L(f(x),f(x'))\le b\le c\|x-x'\|\le c\lambda\|x-x'\|$. Assume
that $x\in\Delta^k\subset K$, $x'\in\Delta^{k'}\subset K$, and
$\|x-x'\|<\sqrt{\frac{2}{n+1}}$ where $\Delta^k,\Delta^{k'}$ are
simplices. Then $\Delta^k\cap\Delta^{k'}\ne\emptyset$.  Hence
there is a geodesic segment $J$ in $\Delta^k\cup\Delta^{k'}$
joining $x$ with $x'$ of length $|J|$. We may assume that it is
piece-wise linear (actually it consists of two straight
intervals). Since $f$ is $\lambda$-Lipschitz on each of the
segments of $J$, we obtain $d_L(f(x),f(x'))\le\lambda |J|\le
\lambda c\|x-x'\|$.\qed
\enddemo

Let $f:X\to\Delta$ be a map of a metric space to the standard
simplex. Denote by
$$
\delta(f)=\inf\{Lip(\psi)\mid \psi:X\to\p\Delta,\
\psi|_{f^{-1}(\p\Delta)}=f|_{f^{-1}(\p\Delta)}\}.
$$

Let $X$ be a proper geodesic metric space with an anti-\v Cech
approximation $\{p_i:X\to N_i\}$ where $\dim N_i=n$ for all $i$.
For every $i$ we define
$$\delta^i=\overline{\lim_{\Delta\to\infty}}\delta(p_i|_{p_i^{-1}(\Delta)})
$$
where $\Delta$ runs over all $n$-simplices from $N_i$.

\proclaim{Proposition 3.2} Suppose that $\as X=n$ for a proper
geodesic metric space $X$ and let $\{p_i:X\to N_i\}$ to be an
anti-\v Cech approximation of $X$ by $n$-dimensional polyhedra.
Then there is $c>0$ such that $\delta^i>c$ for all $i$.
\endproclaim
\demo{Proof} Assume the contrary. By passing to a subsequence we
may assume that $\delta^i\to 0$. Let $C_i\subset N_i$ be a finite
subcomplex such that $\delta(p_i|_{p_i^{-1}(\Delta)})<2\delta^i$
for all $n$-simplices $\Delta\subset N_i\setminus C_i$. Then by
sweeping one can define a map $\xi_i:X\to N_i^{(n-1)}\cup C_i$
which is $2\delta^i$-Lipschitz on every set of the from
$p_i^{-1}(\Delta)$. Since $X$ is geodesic, the map $\xi_i$ is
$2\delta^i$-Lipschitz where $N_i^{(n-1)}\cup C_i$ is taken with
restricted metric from $N_i$ and $N_i$ is supplied with the
uniform geodesic metric. Consider the quotient map
$q_i:N_i^{(n-1)}\cup C_i\to N_i^{(n-1)}\cup C_i/C_i=K_i$ and take
a uniform bounded metric on $K_i\subset \ell_2(K^{(0)})$. By
Proposition 3.1 $q_i$ is $c$-Lipschitz where $c$ depends on $n$
only. Thus, $\psi_i=q_i\circ\xi_i:X\to K_i$ is
$2c\delta^i$-Lipschitz. Therefore, $\{\psi_i:X\to K_i\}$ is an
anti-\v Cech approximation $\{\psi_i:X\to K_i\}$ of $X$ by
$n-1$-dimensional polyhedra. Then by the definition $\as X\le n-1$
(see the discussion after Lemma 2.1), which contradicts to the
assumption of the Proposition.\qed
\enddemo

The following lemma is taken from [Dr4], Lemma 2.2.
\proclaim{Lemma 3.3} Suppose that $X$ and $Y$ are finite uniform
simplicial complexes. Then for every $\lambda$ there exists
$\mu=\mu(\lambda)$ such that every null homotopic
$\lambda$-Lipschitz map $f:X\to Y$ admits a $\mu$-Lipschitz
homotopy $H:X\times I\to Y$ to a constant map.
\endproclaim

Here we give a short review of the classical extension theory. Let
$(K,A)$ be a CW complex pair. Suppose that an extension problem
$$
\CD A @>f>> Y\\
@V{\subset}VV @.\\
K\\
\endCD
$$
for a simply connected space $Y$
is resolved on the $n$-skeleton $K^{(n)}$ of $K$ by a map $g:A\cup
K^{(n)}\to Y$. The map $g$ defines the obstruction cochain
$c_g:C_{n+1}(K,A)\to\pi_n(Y)$ which is a relative cocycle. If a
corresponding relative cohomology class is zero, i.e., if
$c_g=\delta\gamma$ for some $\gamma$, then one can change the map
$g$ on the relative $n$-skeleton without changing it on the
$n-1$-dimensional skeleton such that a new map $g'$ can be
extended to a map $\bar g:A\cup K^{(n+1)}\to Y$ [Hu].

If the homotopy group $\pi_n(Y)$ is supplied with a norm and $Y$
is given some metric one can bring a quantitative statement.

\proclaim{Proposition 3.4} Suppose that $K$ is a simplicial
complex supplied with a uniform metric and let $Y=S^n$ be the unit
$n$-sphere. Suppose that in the above extension problem $g$ is
$\lambda$-Lipschitz and $|\gamma|\le b$. Then $g'$ can be taken
$\mu'$-Lipschitz and $\bar g$ with a $\mu$-Lipschitz restriction
to every $n+1$-simplex where $\mu',\bar\mu$ depend on $n,\lambda$
and $b$ only.
\endproclaim
\demo{Proof} We fix a triangulation on $S^n$, say, by identifying
$S^n$ with the boundary of the standard $n+1$-simplex. Let
$\sigma\subset S^n$ be an $n$-face. Using a simplicial
approximation we may assume that $g$ is simplicial with respect to
some iterated barycentric subdivision $\beta^i(A\cup K^{(n)})$
where $i$ depends on $\lambda$ (and $n$). The Obstruction Theory [Hu]
prescribes a construction of  $g'$ as a map having degree on each
$n$-simplex $\sigma'$ in $K\setminus A$ equal to
$deg(g|_{\sigma'})-\gamma(\sigma')$ where the degree of
$g|_{\sigma'})$ (as well as the degree of $g'|_{\sigma'}$) is
computed for the map of pairs $g|_{\sigma'}:(\sigma',\p\sigma')\to
(S^n,S^n\setminus Int\sigma)$. Since $|deg(g|_{\sigma'})|$ is
bounded by the number of simplices in $\beta^i(\Delta^n)$. Thus,
$|deg(g|_{\sigma'})-\beta(\sigma')|$ is uniformly bounded by a
number that depends only on $n,\lambda$ and $b$. Since every
$\sigma'$ is isometric to the standard $n$-simplex, this degree
can be realized by a $\mu'$-Lipschitz map $g'$ where $\mu'$
depends on $n,\lambda$ and $b$ only. According to Lemma 3.3 the
extension $\bar g$ of $g'$ can be taken to be $\bar\mu$-Lipschitz
on every $n+1$-dimensional simplex of $K$. \qed
\enddemo

Let $\sqcup \Delta^n_i$ be a disjoint union of (oriented)
$n$-simplices. Then the $n$-cochain $\mu$ that take each simplex
to 1 defines a nonzero element  $\mu\in AH_b^n(\sqcup
\Delta^n_i,\sqcup\p \Delta_i^n)$ which we call {\it the
fundamental class}.

The family $\{A_i\}$ of bounded subsets in a metric space $X$ is called
{\it dispersed} if there is a function $s:\R_+\to\R_+$
tending to infinity such that the family $\{ A_i\setminus B_R(x_0)\}$ is
$s(R)$-disjoint.

\proclaim{Lemma 3.5} Suppose that $\as X=n$ for a proper geodesic
metric space $X$ and let $\{p_i:X\to N_i;p^i_{i+1},q^i_{i+1}\}$ be
a regular anti-\v Cech approximation of $X$ by $n$-dimensional
simplicial complexes. Then there is a sequence of $n$-simplices
$\Delta_i\subset N_i$ such that the cohomology homomorphism
induced by $\phi=\sqcup p_i|_{p_i^{-1}(\Delta_i)}$,
$$
\phi^*:AH_b^n(\sqcup\Delta_i,\sqcup\p\Delta_i)@>>> AH^n(\sqcup
p^{-1}_i(\Delta_i),\sqcup p^{-1}_i(\p\Delta_i)),
$$
takes the fundamental class $\mu$ to a nonzero element.
\endproclaim
\demo{Proof} In view of Proposition 3.2 we can take
$\Delta_i\subset N_i$ such that
$\delta(p_i|_{p_i^{-1}(\Delta_i)})>c/2$ and
$\{p_i^{-1}(\Delta_i)\}$ is dispersed.

Denote by $M^k_i=(p^k_i)^{-1}(\Delta_i)$, and
$K^k_i=(p^k_i)^{-1}(\p\Delta_i)$, $k<i$. According to the condition
(2) of a regular anti-\v Cech approximation $M^k_i$ is a subcomplex
of $N_k$. For every function $\kappa:\N\to\N$ with
$\kappa=\kappa(i)<i$ and $\kappa(i)\to \infty$ we define an open
cover $\sU_{\kappa}$ of $\sqcup p^{-1}_i(\Delta_i)$ as follows
$$
\sU_{\kappa}=\{p^{-1}_{\kappa(i)}(Ost(v,M^{\kappa(i)}_i))\mid v\in
(M^{\kappa(i)}_i)^{(0)},\ i\in\N\}$$ where $Ost(v,M)$ denotes the
open star of a vertex $v$ in a complex $M$. Let $\{\sU_i\}$ be the
family of open covers that forms the above anti-\v Cech
approximation. It is easy to check that for any $j$ and for every
subcomplex $L\subset N_j$, $p^{-1}_j(Ost(u,L)=U\cap p^{-1}_j(L)$
where $U\in\sU_j$ and $u$ a vertex in the nerve $N_j$ that
corresponds to $U$. Then
$\sU_{\kappa}=\cup_i\sU_{\kappa(i)}|_{p^{-1}_i(\Delta_i)}$. Then the
nerve of $\sU_{\kappa}$ coincides with $\sqcup_i M^{\kappa(i)}_i$.
Being the composition of simplicial approximations the map $q^k_i$
is a simplicial approximation of $p^k_i$. Thus, we obtain
$q^k_i(M^k_i)=\Delta_i$ for all $k<i$.

Since $\{p_i^{-1}(\Delta_i)\}$ is dispersed, for every $\sU\in
aCov(\sqcup p^{-1}_i(\Delta_i))$ we obtain that the Lebesgue number
of the restrictions tends to infinity:
$L(\sU|_{p^{-1}_i(\Delta_i)})\to \infty$. From here it is easy to
verify that the family of covers
$$\{\sU_{\kappa}\mid \kappa:\N\to\N,\ \kappa=\kappa(i)<i,\ \kappa(i)\to
\infty\}$$ is cofinal in $aCov(\sqcup p^{-1}_i(\Delta_i))$. Note
that $N(\sU_{\kappa(i)}|_{p_i^{-1}(\p\Delta_i)})=K^{\kappa(i)}_i$.

Assume that $\phi^*(\mu)=0$. Then there is $\kappa:\N\to\N$,
$\kappa=\kappa(i)<i$ and $\kappa(i)\to \infty$ such that the
homomorphism $\psi^*$ induced by the simplicial map $$\psi=\sqcup
q^{\kappa(i)}_i|_{M^{\kappa(i)}_i}:\sqcup(M^{\kappa(i)}_i,K^{\kappa(i)}_i)
 \to\sqcup(\Delta_i,\p\Delta_i)$$ takes $\mu$ to zero. Let
$\psi_i=q^{\kappa(i)}_i|_{M_i^{\kappa(i)}}$. Denote by $\mu_i$ the
image of the cocycle $1_i\in C^n(\Delta_i,\p\Delta_i)\cong\Z$ that
takes $\Delta_i$ to 1 under $\psi_i^*$. If $\phi^*_1(\mu)=0$ then
for all sufficiently large $i$ there are a number $b>0$ and
$b$-bounded $n-1$-dimensional cochains $\gamma_i\in
C^{n-1}(M_i^{\kappa(i)})$ such that $\delta \gamma_i=\mu_i$.
Since $q^{k(i)}_i$ is a simplicial approximation of
$p^{\kappa(i)}_i$, we obtain that
$\psi_i(K_i^{\kappa(i)})\subset\p\Delta_i$. We consider the
extension problem:
$$
\CD K_i^{\kappa(i)} @>{\psi_i|_{...}}>>
\p\Delta\\
@VjVV @.\\
M_i^{\kappa(i)}\\
\endCD
$$
where $j:K_i^{\kappa(i)}\to M_i^{\kappa(i)}$ is the inclusion. The
map $\psi_i$ solves this problem on the $n-1$-skeleton and the
relative cocycle $\mu_i$ is the obstruction cocycle to the
solution of the problem on $M_i^{\kappa(i)}$. Since the
obstruction cocycle is a coboundary, by Proposition 3.4 one can
change the map $\psi_i$ on the $(n-1)$-skeleton of
$M_i^{\kappa(i)}\setminus K_i^{\kappa(i)}$ without changing it on
the $(n-2)$-skeleton such that new map has an extension to
$M_i^{\kappa(i)}$. Moreover, we may assume that there is
$\bar\lambda$ such that for every $i$ there is an extension
$\bar\psi_i:M_i^{\kappa(i)}\to\p\Delta_i$ which is
$\bar\lambda$-Lipschitz on every $n$-simplex of $M_i^{\kappa(i)}$.
By Proposition 3.1 we can conclude that $\bar\psi_i$ is
$c'\bar\lambda$-Lipschitz on $M_i^{\kappa(i)}$ for all $i$. Note
that
$$c/2<\delta(\p_i|_{p_i^{-1}(\Delta_i)})\le Lip(\bar\psi_i\circ
p_{\kappa(i)})<c'\bar\lambda\epsilon_{\kappa(i)}.$$ This contradict
to the fact that $\epsilon_{\kappa(i)}\to 0$.\qed
\enddemo

{\bf Asymptotic Alexandroff Theorem.} The classical Alexandroff
Theorem states that $\dim X=\dim_{\Z}X$ for compact metric spaces
provided $\dim X<\infty$. Here we prove the following.
\proclaim{Theorem 3.6} Suppose $\as X<\infty$ for a proper
geodesic metric space $X$. Then $\as_{\Z}X=\as X$.
\endproclaim
\demo{Proof} The inequality $\as_{\Z}X\le\as X$ follows from
Proposition 2.1 and the definition of $\as_{\Z}X$.

Assume that $\as X=n$. Then by Lemma 3.5 $AH^n(Z,Y)\ne 0$ for some
$Y\subset Z\subset X$. Therefore, $\as_{\Z}X\ge n$. \qed
\enddemo

\proclaim{Lemma 3.7} Suppose $\as X<\infty$ for a uniformly path
connected metric space $X$ and let $\as_GX=n$ for some abelian group
$G$. Then there exists a dispersed family $\{U_i\}$ of bounded open
sets in $X$ such that $AH^n(\sqcup\bar U_i,\sqcup\p U_i;G)\ne 0$.
\endproclaim
\demo{Proof} Let $\as X=n$ and let $W\subset X$ be an open subset
such that $AH^n(X,X\setminus W;G)\ne 0$. There are dispersed
families of bounded open sets $\sV^k=\{U^k_i\}$, $k=0,\dots,n$, such
that the union $\sV=\cup_{k=0}^n\sV^k$ is a cover $\sV\in aCov(X)$
(see [DKU] for the construction). Denote by $\sU^k=\sV^k|_W$. Using
induction and Propositions 2.4 and 2.5 we can derive that
$AH^n(X,X\setminus\cup_{U\in\sU^k}U;G)\ne 0$ for some $k$. We take
$\sU^k$ as the desired family $\{U_i\}$. The excision (Proposition
2.6) implies that $AH^n(\sqcup\bar U_i,\sqcup\p U_i;G)\ne 0$. \qed
\enddemo

\head \S4 Relation to the coarse cohomological dimension\endhead

In this section we consider only metric spaces $Y$ that admit an
anti-\v Cech approximation. Using an anti-\v{C}ech approximation
$\{N_i,q^i_{i+1}\}$ of $Y$ John Roe [Ro1] defined coarse cohomology
group $HX^*(Y)$ of $Y$ as the homology of the inverse limit of the
cochain complexes $C^*_0(N_i)$ that consist of cochains with compact
supports. Then the standard argument shows that coarse cohomology
groups fit in the short exact sequence [Ro1]:
$$
0\to {\lim}^1 H_c^{k-1}(N_i)\to HX^k(Y)\to \lim_{\leftarrow}
H_c^k(N_i)\to 0
$$
where $H_c$ stands for the cohomology with compact supports. Let
$A\subset Y$, denote by $A_i=St(p_i(A),N_i)$. We may assume that
$q^i_{i+1}(A_i)\subset A_{i+1}$. Then for relative coarse cohomology
there is an exact sequence:
$$
(*)\ \ 0\to {\lim}^1 H_c^{k-1}(N_i\setminus A_i)\to HX^k(Y,A)\to
\lim_{\leftarrow} H_c^k(N_i\setminus A_i)\to 0.
$$
The coarse cohomology and this exact sequence are defined for
any coefficient group.
Using coarse cohomology with integral coefficients one can define a {\it coarse cohomological
dimension} \cite{Dr1}
$$
\X{Y}=\max\{k\mid HX^k(Y,A)\ne 0\ \text{for some}\ A\subset
Y\}-1.
$$
The shift by one is needed to get the equality $X{-}\dim\R^n=n$.
We note that in \cite{Dr1} this dimension was defined under a
different name. Also in \cite{Dr1} it was suggested to make a
shift in the grading of the coarse cohomology in order to achieve
the equality $X{-}\dim\R^n=n$. Here we embed this shift in the
definition of $X{-}\dim$.

\proclaim{Proposition 4.1} Let $Y$ be a proper metric space, then
$$
\X Y\le\as Y.
$$
\endproclaim
\demo{Proof} Assume that $\as Y\le n$. There is an anti-\v{C}ech
approximation of $Y$ with polyhedra $N_i$ of dimension $\le n$.
Then $H^k_c(N_i,A_i)=0$ for any $A_i$ for $k>n$. Therefore,
$HX^k(Y,A)=0$ for $k>n+1$. Hence $X{-}\dim Y\le n$.\qed
\enddemo
\proclaim{Proposition 4.2} For every metric space $Y$, $$
\X(Y\times\R)=\X Y+1.$$
\endproclaim
\demo{Proof} Let $\X Y=n$ and let $HX^n(Y,A)\ne 0$. Since
for all $k$,
$$H_c^k(N_i\setminus A_i)=H_c^{k+1}((N_i\setminus A_i)\times\R),$$
we obtain that $HX^{n+1}(Y\times R, A\times\R)\ne 0$. Thus,
$\X(Y\times\R)\ge \X Y+1$. In view of the inequality
$\as(Y\times\R)\le \X Y+1$, Proposition 4.1 implies the
inequality $\X(Y\times\R)\le \X Y+1$.
 \qed
\enddemo

\

We note that the coarse cohomological dimension can be defined with
any (abelian) coefficient groups $G$. We will use the notation
${\X}_G Y$.

Let $\{G_i,\phi^{i+1}_i\}$ be an inverse system with bonding maps
$\phi^{i+1}_i:G_{i+1}\to G_i$. For $k>i$ we denote by
$\phi^k_i:G_k\to G_i$ the composition
$\phi_i^{i+1}\circ\dots\circ\phi^k_{k-1}$ and denote by
$$\phi^{\infty}_i:\lim_{\leftarrow}G_j\to G_i$$
the projection from the limit space to the $i$th factor.

\proclaim{Proposition 4.3} Suppose that for an inverse sequence of
countable groups
$${\lim_{\leftarrow}}^1\{G_i,\phi^{i+1}_i\}\ne 0.$$
Then there is $i_0$ such that for every $i\ge i_0$ there is an
element $\alpha_i\in G_i$ with $\phi^i_{i_0}(\alpha_i)\ne 0$ and
$\alpha_i\notin Im(\phi^{\infty}_i)$.
\endproclaim
\demo{Proof} Since $\lim^1G_i$ is nonzero, the Mittag-Lefler
condition does not hold for the system. Therefore there is $i_0$
such that the nested sequence
$$
Im(\phi^{i_0+1}_{i_0})\supset
Im(\phi^{i_0+2}_{i_0})\supset Im(\phi^{i_0+3}_{i_0})\supset\dots
$$
does not stabilize. In view of this we can take $\alpha_i$ such that
$\phi^i_{i_0}(\alpha_i)\in Im(\phi^i_{i_0})\setminus
Im(\phi^k_{i_0})$ for some $k$.  \qed
\enddemo

\proclaim{Proposition 4.4} Suppose that a geodesic metric space
$Y$ that admits an anti-\v Cech approximation $\{N_i\}$ satisfies
the equality $X{-}\dim_G Y=n$ and let $HX^{n+1}(Y,A;G)\ne 0$ for
some $A$ and some countable group $G$. Then the homomorphism
$${\lim_{\leftarrow}}^1H^{n}_c(N_i\setminus A_i;G) @>\cong>> HX^{n+1}(X,A;G)$$
in the short exact sequence (*) is an isomorphism.
\endproclaim
\demo{Proof} We assume that all $N_i$ are given the uniform
geodesic metric.

Assume the contrary. Let
$\alpha\in\lim_{\leftarrow}H^{n+1}_c(N_i\setminus A_i;G)$ and
$\alpha\ne 0$. Let $(\alpha_i)$ be a thread representing $\alpha$,
$\alpha_i\in H^{n+1}_c(N_i\setminus A_i;G)$. From the definition of
cohomology with compact supports it follows that there is a bounded
open set $U_1\subset N_1\setminus A_1$ and an element $\gamma_1\in
H^{n+1}_c(U_1;G)$ which is taken by the inclusion homomorphism to
$\alpha_1$. There is $r_1>0$ such that $U_1$ lies in the
$r_1/2$-neighborhood of $A_1$, $U_1\subset N_{r_1/2}(A_1)$. Let
$W_1=p^{-1}_1(U_1)$.

Since the projection $p_1:X\to N_1$ is uniformly cobounded, the
set $A^1=p_1^{-1}(N_{r_1}(A_1))$ is in bounded distance to $A$.
Then $HX^{n+1}(Y,A^1;G)\cong HX^{n+1}(Y,A;G)$. Moreover for large
enough $k$ there is an isomorphism of inverse sequences
$$
\CD
H^{n+1}_c(N_1\setminus A_k;G) @<<< H^{n+1}(N_1\setminus
A_{k+1};G) @<<< H^{n+1}(N_1\setminus
A_{k+2};G)@<<<\dots\\
@AAA  @AAA @AAA\\
H^{n+1}_c(N_1\setminus A_k^1;G) @<<< H^{n+1}(N_1\setminus
A_{k+1}^1;G) @<<< H^{n+1}(N_1\setminus
A_{k+2}^1;G)@<<<\dots.\\
\endCD
$$
There is a bounded open set $U_2\subset N_k\setminus A_k^1$ and an
element $\gamma_2\in H^{n+1}_c(U_2;G)$ that goes to $\alpha_k$
under the inclusion homomorphism. Let $U_1^2=(q^2_1)^{-1}(U_2)$
and $W_2=p^{-1}_1(U_2^1)$. The commutative diagram
$$
\CD H^{n+1}_c(N_1\setminus A_1;G) @<<< H^{n+1}(N_1\setminus
A_k;G)\\
@AAA @AAA\\
H^{n+1}_c(U^2_1;G) @<(q^2_1|)^*<< H^{n+1}_c(U_2;G)\\
\endCD
$$
implies that $(q^2_1|)^*(\gamma_2)\ne 0$. There is $r_2>0$ such
that $U_2$ lies in $r_2/2$-neighborhood of $A_k$ and we continue
in a similar fashion. As the result we construct a dispersed
sequence of bounded open sets $W_i\subset Y$ and an anti-\v Cech
approximation
$$
(N_1,N_1\setminus(\sqcup_j U^j_1))
@>>>(N_{k_2},N_{k_2}\setminus(\sqcup_j U^j_2))
@>>>(N_{k_3},N_{k_3}\setminus(\sqcup_j U^j_3)) @>>>\dots
$$
of the pair $(Y,Y\setminus(\sqcup_i W_i))$ such that for all $i$
there are nontrivial elements $\gamma_i\in H_c^{n+1}(U_i;G)$ that
survive after translation to the first level. We may assume that
these elements do not belong to the image of the projection from
the higher level. This implies that the system is not
Mittag-Lefler and hence $\lim^1H^{n+1}_c(U_i;G)\ne 0$. Therefore,
$HX^{n+2}(Y,Y\setminus(\sqcup W_i);G)\ne 0$ and hence
${\X}_GY\ge n+1$. This contradicts to the assumption.\qed
\enddemo

A metric space $Y$ is {\it uniformly $n$-connected} if there is a
control function $\rho:\R_+\to\R_+$, $\rho(t)\ge t$, such that for
every $t>0$ and every $y\in Y$ the inclusion $B_t(y)\to
B_{\rho(t)}(y)$ induces zero homomorphism of $k$-dimensional
homotopy groups for $k\le n$. It is called {\it uniformly
contractible} if there is a function $\rho:\R_+\to\R_+$, $\rho(t)\ge
t$, such that for every $t>0$ and every $y\in Y$ the inclusion
$B_t(y)\to B_{\rho(t)}(y)$ is null-homotopic.

We recall that for every $Y$ there is a natural through homomorphism
$c:HX^k(Y;G)\to H_c^k(Y;G)$ [Ro1]:
$$HX^k(Y;G)\to\lim_{\leftarrow} H_c^k(N_i;G)\to H_c^k(Y;G).$$

\proclaim{Theorem 4.5 ([Ro1], page 33)} For a uniformly
$n$-connected metric space $Y$ the map $c:HX^k(Y;G)\to H_c^k(Y;G)$
is an isomorphism for $k\le n$ and for any abelian coefficients
group $G$.
\endproclaim
We denote the global cohomological dimension of a space $Y$ with
respect to the coefficient group $G$ by
$$
gcd_GY=\max\{n\mid H^n_c(Y;G)\ne 0\}.
$$

\proclaim{Proposition 4.6} Let $Y$ be a uniformly $n$-connected
metric space of bounded geometry and let $gld_GY\le n$ for a
countable group $G$. Then ${\X}_GY\ge gld_GY$.
\endproclaim
\demo{Proof} Assume the contrary: ${\X}_GY=m<k=gld_GY$. In view
of Theorem 4.5 $HX^k(Y;G)=H_c^k(Y;G)\ne 0$ and hence, $m\ge k-1$.
Then $m=k-1$. Since $Y$ has bounded geometry, it has an anti-\v Cech
approximation $\{N_i\}$. By Proposition 4.4 and Theorem 4.5 we
obtain
$$
{\lim_{\leftarrow}}^1H^{k}_c(N_i;G)=HX^k(Y;G)=H_c^k(Y;G).
$$
The last group is countable as the direct limit of countable groups.
The first group is uncountable as a nonzero group which is a
$\lim^1$ group of a sequence of countable groups [Ha]. We arrived to
a contradiction.\qed
\enddemo

\proclaim{Proposition 4.7} Suppose that ${\X}_G Y=n$ for a
proper geodesic metric space $Y$ and a countable group $G$. Then
there is a dispersed sequence $\{U_i\}$ of bounded open sets in $Y$
such that $HX^{n+1}(\sqcup \bar U_i,\sqcup\p U_i;G)\ne 0$. Moreover,
there is an anti-\v Cech approximation $\{N_i,p_i\}$ of $Y$ such
that $U_i=p^{-1}_i(V_i)$, $V_i=Int(L_i)$, and $L_i\subset N_i$ is a
finite subcomplex.
\endproclaim
\demo{Proof} Assume that $HX^{n+1}(Y,A;G)\ne 0$. Let
$$
N_1 @ >{q^1_2}>> N_2 @>{q^2_3}>> N_3 @>{q^3_4}>>\dots N_i
@>{q^i_{i+1}}>>\dots
$$
be a regular anti-\v Cech approximation of $Y$. By Proposition 4.4,
$${\lim_{\leftarrow}}^1H^n(N_i,A_i;G)\ne 0.$$ Note that
$q^i_{i+1}:(N_i,A_i)\to(N_{i+1},A_{i+1})$ is homotopic to
$p^i_{i+1}:(N_i,A_i)\to(N_{i+1},A_{i+1})$. Thus,
$(q^i_{i+1})^*=(p^i_{i+1})^*$.

We apply Proposition 4.3 to obtain $i_0$ and $\alpha_i\in
H^n(N_i,A_i;G)$. Without loss of generality we may assume that
$i_0=1$. Then $(p^1_i)^*(\alpha_i)\ne 0$ and $\alpha_i\notin
Im(p^i_{\infty})^*$ where
$$(p^i_{\infty})^*:\lim_{\leftarrow}H^n(N_k,A_k;G)\to
H^n(N_i,A_i;G)$$ is the projection from the limit to $i$th factor in
the inverse sequence $$\{H^n(N_k,A_k;G),(p^k_{k+1})^*\}.$$ Since
$H^n(N_k,A_k;G)=H_c^n(N_k\setminus A_k;G)$ there is a  open set
$V_k\subset N_k\setminus A_k)$ with compact closure $\bar V_k$ and
an element $\beta_k\in H^n_c(V_k;G)$ such that $\beta_k$ goes to
$\alpha_k$ under the inclusion homomorphism. We may assume that
$V=Int L_k$ where $L_k$ is a finite subcomplex of $N_k$. We define
$U_i=p_i^{-1}(V_i)$. For $j\ge i$ we denote by
$L^j_i=(p^i_j)^{-1}(L_j)$. For $j<i$ we denote $L^j_i=N(\sU_j|_{\bar
U_j})$, the nerve of the cover generating $N_j$ restricted to $\bar
U_j$. In view of (*) it suffices to show that $\lim^1H^n(\sqcup_j
(L^j_i,\p L^j_i);G)\ne 0$. The inverse sequence $\{H^n(\sqcup_{j\ge
1}(L^j_i,\p L^j_i);G)\}$ can be mapped epimorphically onto the
sequence $\{H^n_c(\sqcup_{j\ge i}(L^j_i,\p L^j_i);G)\}$. In view of
the 6-term exact sequence for the inverse limit it suffices to show
that $${\lim_{\leftarrow}}^1\{H^n_c(\sqcup_{j\ge i}(L^j_i,\p
L^j_i);G)\}\ne 0.$$ Since for $j>i$ the element $\beta_j$ goes to
nonzero, the result follows.\qed
\enddemo

\

{\bf Coarse cohomological dimension vs asymptotic.} On subgroups of
the reals $\R$ we consider a natural norm $|\ |$. On the mod $p$
group $\Z_p$ we consider the zero semi-norm.

\proclaim{Theorem 4.8} The following holds for every proper metric
space $Y$ with finite asymptotic dimension:\roster
\item{} $\as_GY\ge {\X}_GY$ for every subgroup $G\subset\Q$;
\item{} $\as_{\Z_p}Y={\X}_{\Z_p} Y$ for all $p$.
\endroster
\endproclaim
\demo{Proof} (1) If $G=\Z$, the result follows from Theorem 3.6
and Proposition 4.1.

Now we assume that $G$ is $p$-divisible for some $p$. Let
${\X}_GY=n$ and let $U_i$ be as in Proposition 4.5. Let
$$
\sqcup_i(N_1^i,K_1^i) @>q^1_2>> \sqcup_i(N_2^i,K_2^i) @>q^2_3>>
\sqcup_i(N_3^i,K_3^i) @>>>\dots
$$
be an anti-\v Cech approximation of $\sqcup_i(\bar U_i,\p U_i)$
where $N_i$ are finite complexes. Then
$\lim^1H^n_c(\sqcup_i(N_j^i\setminus K_j^i);G)\ne 0$. Note that each
group $H^n_c(\sqcup_i(N_j^i\setminus K^i_j);G)=\oplus_i
H^n(N_j^i,K_j^i;G)$ is countable. Hence the system is not
Mittag-Lefler. Without loss of generality, we may assume that for
any $m$ there are $i\ge m$ and $k=k(i)$ and an element $\gamma_i\in
H^n(N_{k(i)}^i,K_{k(i)}^i;G)$ such that $(q^1_{k(i)})^*(\gamma_i)\ne
0$ and $\gamma_i\notin Im(q^{k(i)}_{k(i)+1})$. By dividing
$\gamma_i$ by some power of $p$ we may achieve that $\gamma_i$ is
represented by a cocycle with the norm $\le 1$. This defines a
nontrivial element of $AH_b^n(\sqcup_i(N_{k(i)}^i, K_{k(i)}^i);G)$
and of $AH^n(\sqcup \bar U_i,\sqcup\p U_i;G)$. Thus, $\as_GY\ge n$.

(2) The same argument works to show that $\as_{\Z_p}Y\ge
{\X}_{\Z_p}Y$.

Let $\as_{\Z_p}Y=n$. We apply Lemma 3.7 to obtain a dispersed family
of open set $\{U_i\}$ with $AH^n(\sqcup\bar U_i,\sqcup\p
U_i;\Z_p)\ne 0$.  Without loss of generality we may assume that
$U_i=p_i^{-1}(V_i)$ where $\{N_i,p_i\}$ is an anti-\v Cech
approximation and $\bar V_i\subset N_i$ are subcomplexes.
Furthermore, there is $i_0$ and
$$\alpha\in H^n_b( \sqcup_{i\ge i_0}(\bar V_i,\p V_i);\Z_p)=
H^n(\sqcup_{i\ge i_0}(\bar V_i,\p V_i);\Z_p)=\prod_{i\ge i_0}
H^n(\bar V_i,\p V_i;\Z_p)$$ such that $\alpha$ defines a nonzero
element of $AH^n(\sqcup\bar U_i,\sqcup\p U_i;\Z_p)\ne 0$. If
$\alpha=(\alpha_i)_{i\ge i_0}$, this implies that for every function
$\kappa:\N\to \N$, $\kappa(i)>i$, the sequence
$((p_i^{\kappa(i)})^*(\alpha_i))_{i\ge i_0}$ is not eventually zero.
We show that there is $k$ such that for $i>k$ the image
$(p_i^k)^*(\alpha_i)$ is nonzero for infinitely many $i$.  Let
$J_k=\{i\in\N\mid (p_i^k)^*(\alpha_i)\ne 0\}$. Note that $J_k\subset
J_{k+1}$ and $i\in J_i$. If each $J_k$ is finite we can define
$\kappa:\N_{\ge i_0}\to N_{\ge i_0}$ by the formula
$\kappa(J_k\setminus J_{k-1})=k-1$. Then $\kappa\to\infty$ and the
sequence $((p_i^{\kappa(i)})^*(\alpha_i))_{i\ge i_0}$ is zero which
contradicts to the assumption.

Then it follows that
$\lim^1_{\leftarrow}H^n_c(\cup_{j>i}(p^j_i)^{-1}(V_i);\Z_p)\ne 0$.
Therefore $HX^{n+1}(Y, Y\setminus\cup U_i;\Z_p)\ne 0$ and hence
${\X}_{\Z_p}Y\ge n$. Thus, $\as_{\Z_p}Y\le {\X}_{\Z_p}Y$.
 \qed
\enddemo

{\bf Applications to dimension of discrete groups.}

\proclaim{Theorem 4.9} Let $\Gamma$ be a discrete group with finite
$B\Gamma$, then  $\X\Gamma= cd(\Gamma)$.
\endproclaim
\demo{Proof} First we show that $\X\Gamma\ge cd(\Gamma)$. Let
$cd(\Gamma)=n$.  Since $cd(\Gamma)=gld_{\Z}\Gamma$ [Br] and
$E\Gamma$ is uniformly contractible, by Proposition 4.6 we obtain
$\X_{\Z}E\Gamma=\X_{\Z}\Gamma \ge n$.

Let $\X\Gamma= n$ and let $cd\Gamma<n$. By crossing $\Gamma$
with $\Z$ and applying Proposition 4.2 we may assume that $n\ge
3$. Then we have $\dim E\Gamma<n$. Let
$\{N_i,p^i_{i+1},q^i_{i+1}\}$ be a regular anti-\v Cech
approximation of $E\Gamma$. By Proposition 4.7 there is a
dispersed sequence of open sets $U_i=p^{-1}_iV_i$ and elements
$\alpha_i\in H^n_c(V_i)$ with $(p^1_i)^*(\alpha_i)\ne 0$. Let
$W_i$ be the regular neighborhood of $\bar V_i$ in $N_i$, i.e.,
$W_i$ is the star neighborhood of $\bar V_i$ in the second
barycentric subdivision of $N_i$. Let $A_i$ be the regular
neighborhood of $\p V_i$. Since the pair $(W_i,A_i)$ is homotopy
equivalent to $(\bar V_i,\p V_i)$, the element $\alpha_i$ lives in
$H^n(W_i,A_i)$. Since $E\Gamma$ is uniformly contractible there
are lifts $s_i:V^1_i=(p^1_i)^{-1}(\bar V_i)\to E\Gamma_i$ such
that $p_1s_i$ is $r$-closed to the identity where $r$ is the same
for all $i$. Since the Lipschitz constant of $p^1_i$ tends to
zero, for large enough $i$ the maps $p^1_i:(\bar V^1_i,\p
V^1_i)\to (W_i,A_i)$ and $p_is_i:(\bar V^1_i,\p V^1_i)\to
(W_i,A_i)$ are homotopic. Hence
$(p^1_i)^*(\alpha_i)=s_i^*p_i^*(\alpha_i)$. Since $\dim
E\Gamma<n$, $s_i^*$ is zero homomorphism and hence
$s_i^*p_i^*(\alpha_i)=0$. We arrived to a contradiction. \qed
\enddemo
We recall that in the group theoretic language the groups with
finite $B\Gamma$ are called the groups of the type $FL$. Also we
recall that a finitely presented group is called of the type FP if
$B\Gamma$ is dominated by a finite complex.

\proclaim{Proposition 4.10} Let $\Gamma$ be a discrete group of the
type $FP$. Then $\as\Gamma\ge cd(\Gamma)$.
\endproclaim
\demo{Proof} We may assume that $\as\Gamma<\infty$. Let
$cd(\Gamma)=n$. Then $H^n(\Gamma;\Z\Gamma)\ne 0$ [Br]. The condition
$\Gamma\in FP$ is equivalent to the existence of a homotopy
domination $r:K\to B\Gamma$ by a finite complex. We may assume that
$r$ induces an isomorphism of the fundamental groups and
$\pi_i(K)=0$ for $1<i\le n$. Then $H^n(\Gamma;\Z\Gamma)=H_c^n(\tilde
K)$ where $\tilde K$ is the universal cover of $K$. Since the space
$\tilde K$ is uniformly $n$-connected, by Proposition 4.6 we obtain
$\X_{\Z}\tilde K\ge gld_{\Z}\tilde K\ge n$. Since $\tilde K$
coarsely equivalent to $\pi_1(K)=\Gamma$, we obtain by Theorem 3.6
and Proposition 4.8(1) that
$\as\Gamma\ge\as_{\Z}\Gamma=\as_{\Z}\tilde K\ge \X_{\Z}\tilde
K\ge n$. \qed
\enddemo
REMARK. The same argument works for the groups $\Gamma$ with finite
$cd(\Gamma)$ and $B\Gamma$ having finite skeleton in each dimension.

We recall that a group $\Gamma$ is of the type $VFP$ if it admits a
subgroup of finite index of the type $FP$ [Br].

\proclaim{Corollary 4.11} $vcd\Gamma\le\as\Gamma$ for groups
$\Gamma$ of the type VFP.
\endproclaim
Famous Stollings-Swan theorem implies the following.
\proclaim{Corollary 4.12} Every group $\Gamma$ of type VFP with
$\as\Gamma\le 1$ is virtually free.
\endproclaim
This result without the  VFP restriction was proven independently by
Januszkiewicz  and Swiatkowski [JS] and Gentimis [G].

\

\proclaim{Theorem 4.13} Suppose that $\as_{\Z}\Gamma=\as_{\Z_p}\Gamma$
for a group $\Gamma$ of FL type with finite asymptotic dimension and some prime
$p$.
Then $\as\Gamma= cd\Gamma$.
\endproclaim
\demo{Proof} By Theorem
4.8(2) and Theorem 3.6,
$\X_{\Z_p}\Gamma=\as_{\Z}\Gamma=\as\Gamma$. By Theorem 4.9
$cd\Gamma=\X_{\Z}\Gamma\ge \X_{\Z_p}\Gamma=\as\Gamma$.
By Proposition 4.10 $\as\Gamma= cd\Gamma$.
 \qed
\enddemo

\head \S5 A counterexample to the asymptotic analog of Morita's
theorem \endhead

The main idea of the construction can be demonstrated on the
following example. Let $p,q$ be two mutually prime numbers and let
$D$ be a 2-disc with two disjoint discs removed from its interior.
Let $S_1^1$ and $S_2^1$ be the boundaries of the removed discs and
let $S^1$ be the external boundary. We consider a free $\Z_p$
action on $S_1^1$ and a free $\Z_q$ action on $S_2^1$. Let $M$
denote the quotient space which is obtained from $D$ by
factorization of $S^1_1$ and $S^1_2$ to the orbit spaces. Since
the equation $mp+nq=1$ has a solution in integers, there is a
retraction of $M$ to $S^1$. The degree of any such retraction
restricted to the circle $S_1^1/\Z_p$ is $m$ and the degree of the
restriction to $S_2^1/\Z_q$ is $n$ for some $m,n$ satisfying
$mp+nq=1$. Thus, for large $p$ and $q$ these degrees have to be
also large. This allows to construct a uniform complex $M$ that
admits a retraction to the 'boundary' and such that every such
retraction has the Lipschitz constant large.

Now we present the construction. Fix two primes $p>q^2$. For every
natural $k$ we define a 2-dimensional complex $M_k$ as follows.
Let $T_p$ be the mapping cylinder of the degree $p$ map
$z^p:S^1\to S^1$. For a subdivisions of $S^1$ into $p^{i+1}$ and
$p^i$ pieces we fix a mapping cylinder triangulation $T_p^i$ on
$T_p$. Fixing orientation on $S^1$ we may assume that the complex
$T_p^i$ is oriented. We consider the union $M^p_k=T^{k-1}_p\cup
T^{k-2}_p\cup\dots\cup T^1_p$ with identification of the image of
$T^i_p$ with the domain of $T^{i-1}_p$. Let $s\in S^1$ be a base
point. We may assume that it is taken by the map $z^p$ to the base
point (consider $s=1$). Let $M'_k=M^p_k\vee M^q_k$ be the wedge of
these complexes with the base vertices located in the domains. Let
$M$ be a mapping cylinder of the map $\phi:S^1\to S^1\vee S^1$
that collapses two points in $S^1$. We may supply $M$ with a
triangulation having the following properties: (1) every vertex
belongs to at most $p$ edges; (2) the domain $S^1\subset M$ of
$\phi$ has 3 edges; (3) one circle in the target space $S^1\vee
S^1\subset M$ of $\phi$ has $(3p)^k$ edges and the other has
$(3q)^k$ edges. Glue $M$ along the target to $M'_k$ to obtain
$\tilde M_k$. Then we consider a 2-simplex $\Delta$ subdivided in
four 2-simplices by middle points of the edges and delete the
interior of the middle simplex $D=\Delta\setminus Int\sigma$. Glue
$\tilde M_k$ to $\p\sigma$ along the domain of $M$ to obtain
$M_k$. We denote $\p M_k=\p \Delta\subset D$ and will refer to
this set as to the boundary of $M_k$. The end circles of $M_p^k$
and $M_q^k$ we call the boundaries of $p$-hole and $q$-hole in
$M_k$ respectively and denote them by $S^1_p$, $S^1_q$. Note that
they are triangles, i.e., they have three edges. Assume that $M_k$
is given a metric of uniform simplicial complex. We fix a map
$\phi_k:M_k\to\Delta$ which  is simplicial with respect to the
midpoint subdivision of $\Delta$ and such that the original
vertices of $\Delta$ have exactly one preimage each. Clearly, this
map is 1/2-Lipschitz for the uniform metric on $M_k$

We define complexes $M_{k,k-1\dots,k-i}$ for $i=0,\dots,k-1$ by
induction on $i$. Assume that the simplicial complex
$M_{k,k-1\dots,k-i}$ is already constructed. For every 2-simplex
$\Delta\subset M_{k,k-1\dots,k-i}$ we delete its interior and glue
instead a copy of $M_{k,k-1\dots,k-i-1}$ along the boundary $\p
M_{k,k-1\dots,k-i-1}$. This defines a map
$$\phi_{k,\dots,k-i-1}:M_{k,\dots,k-i-1}\to M_{k,\dots,k-i}$$
which is 1/2-Lipschitz. For $j>i$ we denote
$\phi^i_j=\phi_{k,\dots,j}\circ\phi_{k,\dots,j-1}\circ\dots\circ\phi_{k,\dots,i}$.

\proclaim{Proposition 5.1}
 The following holds true: \roster
\item{} For every $k$ there is a retraction $r:M_k\to\p M_k$.
\item{} The Lipschitz constant of every retraction $r$ is greater
than $q^k$. \item{} The Lipschitz constant of every map $f:M_k\to
S^1$ having the restriction $f|_{\p M_k}$ with the nonzero degree
$deg(f|_{\p M_k})\ne 0$ is greater than $q^k$. Moreover,
$Lip(f|_{\p M_k\cup S^1_p\cup S^1_q}) \ge q^k$. \item{} The
Lipschitz constant of every map $f:M_{k,\dots,i}\to S^1$ with the
nonzero degree restriction $deg(f|_{\p M_{k,\dots,i}})\ne 0$ is
greater than $q^i$.
\endroster
Here $S^1$ is given a metric of the boundary of the standard
2-simplex.
\endproclaim
\demo{Proof} (1). Note that $M_k$ is homotopy equivalent to a
2-complex obtained from the wedge $S^1_e\vee S_p^1\vee S_q^1$ by
attaching a 2-cell along the loop $\bar ea^{p^k}b^{q^k}$. Let $n$
and $m$ be such natural numbers that $np^k+mq^k=1$. We consider maps
$r_p:S^1_p\to S^1$ and $r_q:S^1_q\to S^1$ of degree $n$ and $m$
respectively. Then the map $id_{S^1}\cup r_p\cup r_q:S^1_e\vee
S^1_p\vee S^1_q\to S^1$ has an extension $r:M_k\to\p M_k$ since the
attaching map composed with it has the degree $0=-1+np^k+mq^k$.

(2). Note that $Lip(r)\ge Lip(r_q)\ge |m|\ge (p/q)^k>q^k$.

(3). Let $d=deg(f|_{\p M_k})$ and let $d_p=deg(f|_{S^1_p})$, and
$d_q=deg(f|_{S^1_q})$. Then $d_pp^k+d_qq^k=d$. If one of the
coefficients $d_p,d_q$ is zero, then $Lip(f)\ge Lip(f|_{\p M_k})\ge
d\ge q^k$. If both are nonzero, then they are divisible by $d$ and
$d_q/d\ge q^k$ as above. Hence, $Lip(f)\ge q^k$.

(4). By induction on $k-i$. For $i=k$, the result is proven in
(2). For $i<k$ we consider two cases. First we consider the case
when the restriction of $f$ to one of the end circles in
$M_{k,\dots,i}$ has nonzero degree. Then the restriction of $f$ to
a copy of $M_i$ satisfies the conditions of (3). Therefore,
$Lip(f|_{M_i})\ge q^i$.

Next we assume that $f$ has zero degree on every end circle of
$M_{k,\dots,1}$. Then, $f$ defines a map $\tilde f:M_{k,\dots,
i+1}\to \p\Delta^2$ which agrees with $f$ on the 1-dimensional
skeleton of $M_{k,\dots,i+1}^{(1)}\subset M_{k,\dots,i}$. By the
induction assumption $Lip(\tilde f)\ge q^{i+1}$. Since $Lip(\tilde
f)=Lip(\tilde f|_{M_{k,\dots,i+1}^{(1)}})$, we obtain $Lip(f)\ge
q^{i+1}$. \qed
\enddemo
Let $q_k:M_k\to\Delta$ be a simplicial approximation of
$\phi_k:M_k\to\Delta$. Let $n_k\in\N$, denote by
$\xi_k:[0,n_k]\to[0,1]$ the continuous map that collapses
$[0,n_k-1]$ to 0 and maps $[n_k-1,n_k]$ isometrically onto $[0,1]$.
Note that $\xi_k$ is a simplicial approximation of an orientation preserving
homeomorphism $[0, n_k]\to [0,1]$ where $[0,n_k]$ is subdivided into
the intervals of length one.
Let $g_k=q_k\times\xi_k:M_k\times[0,n_k]\to\Delta\times[0,1]$.

\proclaim{Proposition 5.2} There is $\lambda>0$ such that for
every $k$ there is $n_k\in\N$ and a $\lambda$-Lipschitz map
$f_k:M_k\times[0,n_k]\to S^2$ such that $f_k|_{\p
(M_k\times[0,n_k])}=g_k|_{\p (M_k\times[0,n_k])}$ where $\p
(M_k\times[0,n_k])=\p M_k\times[0,n_k]\cup M_k\times\{0,n_k\}$,
$S^2=\p(\Delta^2\times[0,1])$ with the $l_1$-product metric.
\endproclaim
\demo{Proof} Let $m_k=\min\{\|\gamma\|\mid
\delta\gamma=q_k^*(1_{\Delta})\}$ where $1_{\Delta}\in
C^2(\Delta)$ is a simplicial cocycle that takes 1 on $\Delta$.
According to Proposition 5.1(1) $m_k<\infty$. Assume that this
$m_k$ is attained on a cochain $\gamma_k$. We take $n_k=m_k!$ and
consider the extension problem
$$
\CD
\p (M_k\times[0,n_k]) @>g_k|>> \p(\Delta\times [0,1])\\
@V{\subset}VV @.\\
M_k\times[0, n_k].\\
\endCD
$$
We note that the product $M_k\times[0,n_k]$ has a natural
structure of a 3-dimensional cell complex where $[0,n_k]$ is
subdivided into unit intervals, and the map $g_k$ takes its
2-skeleton to $\p(\Delta\times [0,1])$. Then the obstruction
cocycle equals $c=g_k^*(1_{\Delta\times[0,1]})$ where
$1_{\Delta\times[0,1]}$ is the cellular 3-cocycle on
$\Delta\times[0,1]$ that takes this 3-cell to one. We construct a
2-cochain $\beta_k\in C^2(M_k\times[0, n_k])$ with
$\delta\beta_k=c$ and with $|\beta_k|\le 4$.

For every 1-simplex $e$ in $M_k$
with $\gamma_k(e)\ne 0$ we set $\beta_k(e\times
[im,im+1])=sgn(\gamma_k(e))$ for the integer $m=n_k/|\gamma_k(e)|$ and
$i=0,\dots,|\gamma_k(e)|-1$, and set $\beta_k(e\times[s,s+1])=0$
for all other cells. Here $sgn(x)=\cases 1\
if\ x>0\\
-1\ if\ x<0\\
0\ if\ x=0.\\
\endcases$

Note that $\beta_k(e\times[0,n_k])=\gamma_k(e)$.
Indeed,
$$
\beta_k(e\times[0,n_k])=\sum_i\beta_k(e\times
[im,im+1])=|\gamma_k(e)|sgn(\gamma_k(e))=\gamma_k(e).$$
For every 2-simplex $\sigma\subset M_i$ we set
$\beta_k(\sigma\times\{l\})=-\beta_k(\p\sigma\times[0,l])$ for all
natural $l< n_k$ and we define
$\beta_k(\sigma\times\{0\})=\beta_k(\sigma\times\{n_k\})=0$.

Show that $\delta\beta_k=c=g_k^*(1_{\Delta\times[0,1]})$. Indeed,
$$\delta\beta_k(\sigma\times[s-1,s])=\beta_k(\p\sigma\times[s-1,s]+
\sigma\times\{s\}-\sigma\times\{s-1\})
=\beta_k(\p\sigma\times[s-1,s])-$$ $$\beta_k(\p\sigma\times[0,s])+
\beta_k(\p\sigma\times[0,s-1])=0=g_k^*(1_{\Delta\times[0,1]})(\sigma\times[s-1,s])$$
for $s<n_k$. For $s=n_k$ we obtain
$$\delta\beta_k(\sigma\times[s-1,s])=\beta_k(\p\sigma\times[s-1,s])+
\beta_k(\p\sigma\times[0,s-1])=
\beta_k(\p\sigma\times[0,n_k])=\gamma_k(\p\sigma)=$$
$$\delta\gamma_k(\sigma)=q_k^*(1_{\Delta})(\sigma)=
(q_k\times\xi)^*(1_{\Delta\times[0,1]})(\sigma\times[s-1,s]).$$

We show that the cochain $\beta_k$ is bounded.
Let
$\p\sigma=a+b+c$ where $a,b,c$ are sides (with signs) of a 2-simplex
$\sigma\subset M_k$.
Note that
$$|\beta_k(\sigma\times \{l\})|=|\beta_k(\p\sigma\times[0,l])|=|\beta_k(a\times[0,l])+
\beta_k(b\times[0,l])+\beta_k(c\times[0,l])|=$$
$$|sgn(\gamma_k(a))[\frac{l|\gamma_k(a)|}{n_k}]+
sign(\gamma_k(b))[\frac{l|\gamma_k(b)|}{n_k}]+
sgn(\gamma_k(c))[\frac{l|\gamma_k(c)|}{n_k}]|\le$$
$$\frac{l}{n_k}(\gamma_k(a)+\gamma_k(b)+\gamma_k(c))+3=
\le \frac{l}{n_k}(q_k^*(1_{\Delta})(\sigma))+3\le 4$$ where $[x]$
denote the integral part of $x$.

By the Obstruction Theory (see Proposition 3.4) the map $g_k$
restricted to the 2-skeleton of $M_k\times[0,n_k]$ can be changed
on 2-cells lying in
$M_k\times[0,n_k]\setminus\p(M_k\times[0,n_k])$ without changing
on the 1-skeleton by means of the cochain $\beta_k$ in such a way
that a new map is $\mu(|\beta_k|)$-Lipschitz and it is extendible
over $M_k\times[0,n_k]$ to a $\lambda$-Lipschitz map
$f_k:M_k\times[0,n_k]\to S^2$.  \qed
\enddemo

\proclaim{Theorem 5.3} There is a proper geodesic metric space $Y$
of bounded geometry with dimensions $\as Y=2$ and $\as
(Y\times\R)=2$.
\endproclaim
\demo{Proof} Let $Y_0$ be the disjoint union of the uniform
complexes $M_{k,\dots,1}$. To make it into a geodesic metric space
we attach $M_{k,\dots,1}$ to the half-line $\R_+$ at the point
$2^k\in\R_+$ for all $k$:
$$
Y=\R_+\cup_kM_{k,\dots,1}.
$$

In view of Theorem 3.6 it suffices to prove that $\as Y\ge 2$. We
show that $$AH^2(\sqcup_k M_{k,\dots,1},\sqcup_k\p
M_{k,\dots,1})\ne 0.$$

We note that for every $k$ the projection $\phi^i_k$ takes
$M_{k,\dots,i}$ to $\Delta^2$ in such a way that it is a
homeomorphism on the boundary $\p M_{k,\dots,i}\cong \p\Delta^2$.
Let $i:\N_+\to\N_+$, $i(k)\le k$, $\lim i(k)=\infty$ and let
$q_i:\sqcup_kM_{k,\dots,i(k)}\to\sqcup_k\Delta^2$ be a simplicial
approximation of $\sqcup_k\phi^i(k)_k$. We show that for every
$i:\N_+\to\N_+$, $i(k)\le k$, $\lim i(k)=\infty$,
$$
q_i^*:AH_b^2(\sqcup_k\Delta^2,\sqcup_k\p\Delta^2)\to
AH_b^2(\sqcup_kM_{k,\dots,i(k)},\sqcup_k\p M_{k,\dots,i(k)})$$
takes the fundamental class $\mu$ to nonzero element. Then
$$
AH_b^2(\sqcup_k\Delta^2,\sqcup_k\p\Delta^2)\to
AH_b^2(\sqcup_kM_{k,\dots,1},\sqcup_k\p M_{k,\dots,1})$$ will be a
nonzero homomorphism.

Indeed, the image $q_i^*(1_{\Delta^2})$ of the fundamental cocycle
is the obstruction cocycle $C_f$ for a retraction
$f:M_{k,\dots,i(k)}\to\p \Delta^2$ defined on the 1-skeleton as
$q_i$. By Proposition 5.1(1) and the Obstruction Theory,
$C_f=\delta\gamma_k$ for some simplicial cochain $\gamma_k\in
C^1(M_{k,\dots,i(k)}$. By Proposition 5.1 (4) and the Obstruction
Theory, $\|\gamma_k\|\ge q^i(k)-1$. Since $i(k)\to\infty$, the
cocycle $q^*(1_{\Delta^2})$ defines a nontrivial element in
$AH_b^2(\sqcup_kM_{k,\dots,i(k)},\sqcup_k\p M_{k,\dots,i(k)})$.

\

Now we show that $\as(Y_0\times\R)\le 2$.  We construct an anti-\v
Cech approximation of $Y_0\times\R$ by 2-dimensional complexes. We
recall that the map $\phi^i_j:M_{k,\dots,j}\to M_{k,\dots, i}$ is
defined for $k\ge j\ge i$. Let $Y_i=M_i\sqcup M_{i+1,i}\sqcup
M_{i+2,i+1,i}\sqcup\dots$. We denote
$$
\psi^1_i=\sqcup_{k\ge i}\phi^1_i:\sqcup_{k\ge i}
M_{k,\dots,1}=Y_0\to Y_i
$$
and
$$\psi_i=\phi_i\sqcup\phi_{i+1,i}\sqcup\phi_{i+2,i+1,i}\sqcup\dots:Y_i\to\Delta\sqcup Y_{i+1}  .$$

We consider the map $\psi_i\times\frac{1}{n_i}:Y_i\times\R\to
(\Delta\sqcup Y_{i+1})\times\R$ where $n_i$ is from Proposition
5.2. We subdivide  $\R$ in the first product into the interval of
the form $[ln_i,l(l+1)n_i]$ and subdivide $\R$ in the second
product into the unit intervals to turn
$\psi_i\times\frac{1}{n_i}$ into a cellular map. By Proposition
5.2 there is a $\lambda$-Lipschitz sweeping
$\xi_i:Y_i\times\R\to((\Delta\sqcup Y_{i+1})\times\R)^{(2)}=N_i$
onto the 2-skeleton. Then the composition
$\xi_i\circ(\psi^1_i\times\frac{1}{i}):Y_0\times\R\to N_i$ is a
uniformly cobounded $\lambda/i$-Lipschitz map onto a 2-dimensional
complex. To make it into a genuine  anti-\v Cech approximation one
needs to triangulate all the prisms in $N_i$. \qed
\enddemo
Theorem 5.3 together with Proposition 4.2 implies
\proclaim{Corollary 5.4} If $Y$ as above, then $X{-}\dim Y=1$.
\endproclaim

\head \S6 A construction of large simplicial complexes \endhead

We present here a construction of infinite locally finite uniform
simplicial complexes. Let
$\{\phi_k:M_k\to\Delta^n,\chi_k:M_k\to\Delta^n\}$ be a sequence of
maps of $n$-dimensional simplicial complexes such that for each $k$,
$\chi_k$ is a light simplicial map (and hence a retraction) and
$\phi_k$ is a simplicial map to a some subdivision $\tau_k$ of the
of $\Delta^n$ such that $mesh(\tau_k)<\delta<1$ for fixed $\delta$
for all $k$. We recall that a map $\phi$ is called {\it light} if
all point preimages are 0-dimensional. A light simplicial map
$\phi:M\to\Delta^n$ is an isomorphism on every $n$-simplex
$\sigma\subset M$. Williams calls the complexes that admit a light
simplicial map onto the simplex $\Delta^n$ as {\it complexes over
$\Delta^n$} [Wi].

Let $s_k:\Delta^n\to M_k$ be sections of $\chi_k$.

Consider the following diagram:
$$
\CD M_1 @<{\to}<<    M_{2,1} @<{\to}<< M_{3,2,1} @<{\to}<<\dots
M_{k,\dots,1} @<{\to}<< M_{k+1,\dots,1} @<{\to}<<\dots\\
@VV\phi_1V  @ VV\phi_{2,1}V  @VV{\phi_{3,2,1}}V   @
V\phi_{k,\dots,1}VV
@V\phi_{k+1,\dots,1}VV\\
\Delta^n @<{\to}<\chi_2<   M_2 @<{\to}<<  M_{3,2} @<{\to}<< \dots
M_{k,\dots,2} @<{\to}<<
M_{k+1,\dots,2}@<{\to}<<\dots\\
@. @V\phi_2VV @V\phi_{3,2}VV @V\phi_{k,\dots,2}VV @V\phi_{k+1,\dots,2}VV\\
@. \Delta^n @<{\to}<\chi_3< M_3\ @<{\to}<<\dots M_{k,\dots,3}
@<{\to}<<
M_{k+1,\dots,3} @<{\to}<<\dots\\
@. @. @. \dots\dots\dots\dots\\
@. @. @. @V\phi_kVV @V\phi_{k,k-1}VV\\
@. @. @. \Delta^n @<{\to}<\chi_{k+1}< M_{k+1} @<{\to}<<\dots\\
@. @. @. @. @V\phi_{k+1}VV\\
@. @. @. @. \Delta^n @<{\to}<<\dots\\
\endCD
$$
where all horizontal left arrows  are generated by the maps $\chi_k$
by taking pull-back, all horizontal right arrows are generated by
sections $s_k:\Delta^n\to M_k$ of maps $\chi_k$, and all squares are
pull-back diagrams. We define simplicial complexes $Y_k$ as the
direct limit of embeddings:
$$
Y_k=\lim_{\rightarrow}\{M_k\to M_{k+1,k}\to M_{k+2,k+1,k}\to\dots\}.
$$
Note that for $k<l$ there are the natural projections
$\phi^k_l:Y_k\to Y_l$ obtained as the direct limits in the above
diagram. We denote the top horizontal sections by
$\sigma_k:M_{k,\dots,1}\to M_{k+1,\dots,1}$. Denote by
$$
Y=Y_1=\lim_{\rightarrow}\{M_1 @>\sigma_1>> M_{2,1} @>\sigma_2>>
M_{3,2,1}\to\dots\}.
$$

We fix a simplicial approximation $\rho_k:\tau_k\to\Delta^n$ of the
identity map $id:|\tau_k|\to\Delta$. It defines a simplicial
approximation $q^{k-1}_k$ of the map $\phi^{k-1}_k$.

We assume that all $Y_k$ are given the uniform geodesic metrics.
\proclaim{Proposition 6.1} The family
$$\{\phi^1_k:Y\to Y_k, q^k_{k+1},\phi^k_{k+1}\}$$ is
a regular anti-\v{C}ech approximation of $Y$.
\endproclaim
\demo{Proof} We consider $\sU_k=(\phi^1_k)^{-1}\{Ost(v, Y_k)\}$. By
the construction $\phi^1_k$ is a simplicial map to a $k-1$-iterated
$\delta$-subdivision of a uniform complex. Therefore, $\phi^1_k$ is
$\delta^{k-1}$-Lipschitz on every simplex of $Y$. Since the metric
on $Y$ is geodesic, it is $\delta^{k-1}$-Lipschitz. Since
$\delta^k\to 0$, the Lebesgue numbers $L(\sU_k)$ tends to infinity.
Since $\phi^k_{k+1}$ is simplicial with respect to a subdivision,
for every vertex $v\in Y_k$ there is a vertex $u\in Y_{k+1}$ such
that
$$\phi^k_{k+1}(Ost(v,Y_k))\subset Ost(\phi^k_{k+1}(v),\tau_{k+1})\subset
Ost(u, Y_{k+1}).$$ Therefore, $Ost(v, Y_k)\subset
(\phi^k_{k+1})^{-1}(Ost(u, Y_{k+1})$ and hence
$$(\phi^1_k)^{-1}(Ost(v, Y_k)\subset (\phi^1_{k+1})^{-1}(Ost(u, Y_{k+1})$$
and the condition $\sU_k\prec\sU_{k+1}$ is checked. \qed
\enddemo
REMARK. The space $Y$ has natural compactification
$$\bar Y=\lim_{\leftarrow}\{M_1\leftarrow M_{2,1}\leftarrow
M_{3,2,1}\leftarrow\dots\}.$$

The complexes $M_k$ constructed in \S5 can be considered to be
oriented. Let $\chi:M_k\to\Delta^2$ be a simplicial map defined by
the orientation. We consider the following triangulation $\tau$ of
the 2-simplex $\Delta^2$: First, we take the midpoint subdivision
of $\Delta^2$ and then take the cone subdivision of the central
2-simplex. Let $\phi_k:(M_k,S^1)\to(\Delta^2,\partial\Delta^2)$ be
a simplicial map to the subdivision $\tau$ of $\Delta^2$ which
takes $M'_k$ to the center  and maps $D$ to $D$ by the identity
map. Then Theorem 5.3 can be stated as follows. \proclaim{Theorem
6.2} The family $\{\phi_k:M_k\to\Delta^2,\chi_k:M_k\to\Delta^2\}$
for any choice of sections $s_k:\Delta^2\to M_k$ defines the space
$Y$ as above with an anti-\v Cech approximation $\{\phi^1_k:Y\to
Y_k\}$ such that $\as (Y\times\R)=\as Y=2$.
\endproclaim
Note that if $s_k(\Delta^2)\cap\p M_k=\emptyset$ for all $k$ then
the boundary of the complex $Y$ is homeomorphic to the circle
$S^1$.

\

\enddocument